\documentstyle[12pt,twoside]{article}

\begin{document}
\newcommand{\myappendix}{\setcounter{equation}{0}\addtocounter{section}{1}\section*{Appendix\thesection}}
\newcommand{\mysection}{\setcounter{equation}{0}\section}
\renewcommand{\theequation}{\thesection.\arabic{equation}}

\thispagestyle{empty}
\topmargin=-.5in
\oddsidemargin=.1in
\evensidemargin=.1in
\vsize=23.5cm
\hsize=16cm
\textheight=23.5cm
\textwidth=16cm

\def\hook#1{{\vrule height#1pt width0.4pt depth0pt}}
\def\leftrighthookfill#1{$\mathsurround=0pt \mathord\hook#1
       \hrulefill\mathord\hook#1$}
\def\underhook#1{\vtop{\ialign{##\crcr                 % |_| under
       $\hfil\displaystyle{#1}\hfil$\crcr
       \noalign{\kern-1pt\nointerlineskip\vskip2pt}
       \leftrighthookfill5\crcr}}}

%\mathbold{e^+e^- \rightarrow f}$\bar{\mathbold{f}}${\bf .}

\font\blackboard=msbm10 at 12pt
\def\zz{\hbox{\blackboard Z}}

\def\be{\begin{equation}}
\def\ee{\end{equation}}
\def\bea{\begin{eqnarray}}
\def\ena{\end{eqnarray}}

\def\nn{\nonumber\\}
\def\up{\uparrow}
\def\da{\downarrow}
\def\imp{~~\Rightarrow~~}
\def\q{\qquad\qquad\qquad}
\def\to{\rightarrow}

\def\frac#1#2{{\textstyle{#1\over#2}}}	   % fraction

\def\ihalf{{\textstyle{i \over 2}}}
\def\half{\frac{1}{2}}
\def\bD{{\bar D}}
\def\bj{{\bar\psi}}
\def\pa{\partial}                              % curly d
\def\de{\nabla}                                       % del
\def\dell{\bigtriangledown} % hi ho the dairy-o
\def\su{\sum}                                         % summation
\def\pr{\prod}                                        % product
\def\mn{{\m\n}}

\def\({\left(} \def\){\right)} \def\<{\langle } \def\>{\rangle }
\def\[{\left[} \def\]{\right]} \def\lb{\left\{} \def\rb{\right\}}

\def\a{\alpha}
\def\b{\beta}
\def\c{\chi}
\def\d{\delta}
\def\e{\in}           % Also, \varepsilon
\def\f{\phi}               %      \varphi
\def\g{\gamma}
\def\h{\eta}
\def\i{\iota}
\def\j{\psi}
\def\k{\kappa}                    % Also, \varkappa (see below)
\def\l{\lambda}
\def\m{\mu}
\def\n{\nu}
\def\o{\omega}
\def\p{\pi}                % Also, \varpi
\def\q{\theta}                    %     \vartheta
\def\r{\rho}                      %     \varrho
\def\s{\sigma}                    %     \varsigma
\def\t{\tau}
\def\u{\upsilon}
\def\x{\xi}
\def\z{\zeta}
\def\D{\Delta}
\def\F{\Phi}
\def\G{\Gamma}
\def\J{\Psi}
\def\L{\Lambda}
\def\O{\Omega}
\def\P{\Pi}
\def\Q{\Theta}
\def\S{\Sigma}
\def\U{\Upsilon}
\def\X{\Xi}
\def\del{\partial}

\def\ca{{\cal A}}
\def\cb{{\cal B}}
\def\cc{{\cal C}}
\def\cd{{\cal D}}
\def\ce{{\cal E}}
\def\cf{{\cal F}}
\def\cg{{\cal G}}
\def\ch{{\cal H}}
\def\ci{{\cal I}}
\def\cj{{\cal J}}
\def\ck{{\cal K}}
\def\cl{{\cal L}}
\def\cm{{\cal M}}
\def\cn{{\cal N}}
\def\co{{\cal O}}
\def\cp{{\cal P}}
\def\cq{{\cal Q}}
\def\car{{\cal R}}
\def\cs{{\cal S}}
\def\ct{{\cal T}}
\def\cu{{\cal U}}
\def\cv{{\cal V}}
\def\cw{{\cal W}}
\def\cx{{\cal X}}
\def\cy{{\cal Y}}
\def\cz{{\cal Z}}

\def\under#1#2{\mathop{\null#2}\limits_{#1}}
\def\iff{\leftrightarrow}                      % <-->
\def\conj{{\hbox{\large *}}} % complex conjugate
\def\lconj{{\hbox{\footnotesize *}}}          % little "
\def\dg{\sp\dagger} % hermitian conjugate
\def\ddg{\sp\ddagger} % double dagger

\def\NP{Nucl.\ Phys.\ B}
\def\PL{Phys.\ Lett. }
\def\PR{Phys.\ Rev.\ Lett. }
\def\PRD{Phys.\ Rev.\ D}

\def\ITP{\footnote{Work supported by National Science Foundation
  grant PHY 97-22101.}\\[.1in] {\it Institute for Theoretical Physics\\
  State University of New York, Stony Brook, NY 11794-3840}\\[.1in]}

%\footnote{\rightskip=-65pt
\newcommand\noi{\noindent}
\newcommand\seq{\;\;=\;\;}

\newcommand\ie {{\it i.e.}}
\newcommand\eg {{\it e.g.}}
\newcommand\etc{{\it etc. }}
\newcommand\viz{{\it viz. }}

\newcommand{\bfl}{\begin{flushleft}}
\newcommand{\efl}{\end{flushleft}}
\newcommand{\bc}{\begin{center}}
\newcommand{\ec}{\end{center}}

\def\fig #1 by #2 (#3){\vbox to #2{
	\hrule width #1 height 0pt depth 0pt\vfill\special{picture #3}}}

\def\boxit#1{\leavevmode\thinspace\hbox{\vrule\vtop{\vbox{
	\hrule\kern1pt\hbox{\vphantom{\bf/}\thinspace{\bf#1}\thinspace}}
	\kern1pt\hrule}\vrule}\thinspace}

% $$ 1{\buildrel<\over{\widetilde{\phantom{m}}}}2 \quad\quad
% 1\mathop<\limits_{\widetilde{\phantom{m}}}2 $$

\def\llsim{\:\lower1.8ex\hbox{$\buildrel<\over{\widetilde{\phantom{
<}}}$}\:}
\def\lsim{\;\raise.4ex\hbox{$\mathop<\limits_{\widetilde{\phantom{<}}}$}\;}
%$$ 1\llsim 2\lsim 3 $$
%$$ 1\le 2\le 3 $$

%\topmargin -.5in
%\oddsidemargin=.1in
%\evensidemargin=.1in
%\vsize=23.5cm
%\hsize=16cm
%\textheight=23.5cm
%\textwidth=16cm
%\baselineskip=16pt
%\pagestyle{empty}

\hfill{ITP-SB-98-49}
\vspace{.50in}

\begin{center}
{\Large \bf The perturbation $\varphi_{2,1}$ of the $M(p, p+1)$ models of
conformal field
theory and related polynomial - character identities}

\vspace{.50in}

Dedicated to George E. Andrews on his 60th birthday.

\vspace{.50in}

{\large by} \\

\vspace{.50in}
Alexander Berkovich\footnote{e-mail: alexb@insti.physics.sunysb.edu}and
Barry M. McCoy\footnote{e-mail: mccoy@max.physics.sunysb.edu}

\vspace{.50in}

{\it Institute for Theoretical Physics\\
State University of New York \\
Stony Brook, NY 11794-3840}

\end{center}

\vspace{.50in}
\centerline{\large Abstract}
\bigskip
Using $q$-trinomial coefficients of Andrews and Baxter along with the
technique of telescopic expansions,
we propose and prove a complete set of polynomial identities of
Rogers-Ramanujan type for $M(p, p+1)$
models of conformal field theory perturbed by the operator $\varphi_{2,1}$.
The bosonic form of our
polynomials is closely related to corner transfer matrix sums which arise
in the computation of the order
parameter in the regime $1^+$ of $A_{p-1}$ dilute models.  In the limit where
the degree of the polynomials tends to
infinity our identities provide new companion fermionic representations for
all Virasoro characters of unitary minimal series.

\eject

\mysection{Introduction}

In 1987, G. Andrews and R. Baxter \cite{bone} introduced 
the very fruitful notion of
the $q$ - trinomial coefficients
\bea
\( \begin{array}{c} L, q \\ A \end{array} \)^n = \left\{  \begin{array}{ll}
\sum_{j \geq  0} q^{j(j+A-n)}
\displaystyle{(q)_L \over (q)_j (q)_{j+A} (q)_{L-2j-A}} \; {\rm for} \;
\mid A \mid \leq L\\ 0 - \; {\rm
otherwise}
\end{array} \right.
\label{oneone}
\ena
where
\bea
(a,q)_k = (a)_k = \left\{ \begin{array}{cl} \displaystyle\prod^{k-1}_{j=0}
(1-aq^j) & {\rm for} \; k
\in \zz_{> 0}, \\ 1 & {\rm for} \; k =0,  \\ \\ {1 \over
\displaystyle\prod^{-k}_{j=1} (1-a
q^{-j})} & {\rm for} \; k \in \zz_{< 0} \end{array} \right.
\label{onetwo}
\ena
and $L \in$ \zz$_{\geq 0}$; $A, n \in$ \zz.\\
The rich properties of these objects were extensively studied in 
\cite{bone}-\cite{bfour}.
Recently, $q$ - trinomial
coefficients were used to generalize classical results of Bailey
\cite{bfive}.  This
generalization \cite{bsix}, termed
``Trinomial Analogue of Bailey's lemma", led these authors and P. Pearce
\cite{bseven} to the following $q$ - series
identities of Rogers-Ramanujan type (which hold for $p \geq 4, p \in$ \zz)
\bea
&& \sum_{\scriptstyle {\bf m}\atop\scriptstyle m_{p-2} \equiv 0({\rm mod}
\: 2)} q^{{1 \over 4} \; {\bf m}
{\tilde C}_{p-1} \; {\bf m}} {1 \over (q)_{m_{0}}} \prod^{p-2}_{j=1} \[
\begin{array}{c} ({1 \over 2}
{\tilde I}_{p-1}{\bf m})_j \\ m_j \end{array} \]_q = \nn
&& \hspace{.50in} = {1 \over (q)_\infty} \sum^\infty_{j=-\infty} \left\{
q^{p(p+1)j^2+j} -
q^{(pj+1)((p+1)j+1)}
\right\}
\label{onethree}
\ena
where ${\tilde I}_{p-1}$ is the matrix of dimension $(p-1) \times (p-1)$:
\bea
&& ({\tilde I}_{p-1})_{0,j} =- ({\tilde I}_{p-1})_{j, 0} =- \delta_{j,2} \;
{\rm for} \;  0 \leq j
\leq p-2\nn
&& ({\tilde I}_{p-1})_{i,j} = \delta_{i,j+1} + \delta_{i,j-1} \; {\rm for}
\; 1 \leq i, j \leq
p-2
\label{onefour}
\ena
where
\bea
\delta_{a,b} = \left\{ \begin{array}{l} 1 \; {\rm if} \; a=b \\ 0 \;\;
{\rm otherwise} \end{array}
\right.
\label{onefive}
\ena
${\tilde C}_{p-1} =2 - {\tilde I}_{p-1}$ and $q$ - binomial coefficients $
\[ \begin{array}{c} n+m \\ n
\end{array} \]_q$ are defined in standard fashion as
\bea
\[ \begin{array}{c} n+m \\ n \end{array} \]_q = \[ \begin{array}{c}  n+m \\
m \end{array} \]_q
= \left\{ \begin{array}{cc} {(q)_{n+m} \over (q)_n (q)_m} & {\rm for} \;
n,m \in \zz_{\geq 0} \\ 0 &
{\rm otherwise} \; .\end{array} \; \right.
\label{onesix}
\ena
Throughout this paper we impose the notation that components of a
vector ${\bf a}$ are either
denoted by $({\bf a})_i$ or $a_i$.  By ${\bf a}_i$ we denote a vector labeled
by $i$ and not its $i^{\rm th}$
component.\\
Recalling the Rocha-Caridy (bosonic) form for the Virasoro characters of $M(p,
p+1)$ minimal models of
conformal field theory \cite{beight} 
\bea
\chi^{p,p+1}_{r,s} (q) = {1 \over (q)_\infty} \sum^\infty_{j=- \infty}
\left\{ q^{pp' j^2 +j
(p'r-ps)} -q^{(pj+r)(p'j+s)}
\right\}
\label{oneseven}
\ena
with $1 \leq r \leq p-1; 1 \leq s \leq p$ and $p'=p+1$, one can immediately
recognize the right-hand-side
of (\ref{onethree}) as $\chi^{p,p+1}_{1,1} (q)$.  Following
well-established tradition \cite{bnine}-\cite{bten}
 we refer to the left-hand-side of
(\ref{onethree}) as a fermionic
representation for
$\chi^{p,p+1}_{1,1} (q)$.\\
It is important to notice that this fermionic representation is quite
different from the one previously
proposed in \cite{bnine} and proven in \cite{beleven}-\cite{belevena}:
\bea
\sum_{\scriptstyle {\bf m}\atop\scriptstyle m_{p-2} \equiv 0({\rm mod} \:
2)} q^{\frac14 {\bf m} \;
C_{p-2} \; {\bf m}} {1 \over (q)_{m{_1}}} \prod^{p-2}_{j=2} \[
\begin{array}{c} ({1 \over 2} I_{p-2} {\bf
m})_j \\ m_j
\end{array} \]_q = \chi^{p,p+1}_{1,1} (q)
\label{oneeight}
\ena
where $C_{p-2} = 2-I_{p-2}$ and $ I_{p-2}$ are Cartan and incidence matrices,
respectively, of the simple Lie
algebras $A_{p-2}$, i.e. \\
\bea
(I_{p-2})_{a,b} = \delta_{a,b+1} + \delta_{a,b-1} \quad {\rm for} \quad 1
\leq a,b \leq p-2 \; .
\label{onenine}
\ena
To understand the meaning of this new representation let us recall, as
was shown by A.
Zamolodchikov \cite{btwelve}, that there exist certain
operators $\{ \varphi_{1,3}; \varphi_{1,2};
\varphi_{2,1} \}$ with which one can perturb $M(p,p+1)$ - models so that
these models will remain exactly
integrable even though they will
lose conformal invariance.  Exact integrability is a very powerful
constraint which often enables one to
solve model completely in terms of so-called thermodynamic Bethe Ansatz
equations \cite{bthirteen}.  We remark that
the incidence matrices such as $\tilde{I}_{p-1}, I_{p-2}$ carry important
information about the structure
of these equations.

One way to recognize particular integrable perturbation from $q$ - series
identities such as
(\ref{onethree}), (\ref{oneeight}) is to find appropriate polynomial
identities which reduce to
the original ones when the degree of polynomials tends to infinity.  This
procedure
was termed ``finitization" in
\cite{bfourteen}.  For $L \in \zz_{\geq 0}$ the finitized version of (\ref{onethree})
connected in \cite{bseven} is
\bea
&& \sum_{\scriptstyle {\bf m}\atop\scriptstyle m_{p-2} \equiv 0({\rm mod}
\: 2)} q^{{1 \over 4} \; {\bf m}
{\tilde C}_{p-1} \; {\bf m}}  \prod^{p-2}_{j=0} \[ \begin{array}{c} ({1
\over 2}
{\tilde I}_{p-1}{\bf m} + L  {\bf e}_0)_j \\ m_j \end{array} \]_q = \nn
&& = \sum^{\infty}_{j=-\infty} \left\{ q^{pp'j^2+j} \( \begin{array}{c} L,
q \\ 2pj \end{array} \)^0
-q^{(jp+1)(jp'+1)} \( \begin{array}{c} L,q \\ 2 pj+2 \end{array} \)^0 \right\}
\label{oneten}
\ena
where the $p-1$ dimensional unit vector in the $a$ - direction ${\bf e}_a$ is
\bea
{\bf e}_a  = \left\{ \begin{array}{l} \delta_{a,b} \; {\rm for} \; a,b
=0,1,2 \ldots, p-2 \\ 0 \quad
{\rm otherwise} \; .\end{array} \right.
\label{oneeleven}
\ena
Bosonic polynomials on the right-hand-side of (\ref{oneten}) emerged
previously in the order parameter
calculations \cite{bfifteen} for regime $1^+$ of dilute $A_{p-1}$ - model.
  It was shown in \cite{bfifteen}
that the scaling limit of regime
 $1^+ A_{p-1}$ - model represents the $M(p,p+1)$ - model perturbed by
$\varphi_{2,1}$.  Thus, the fermionic
representation (\ref{onethree}) for the $\chi^{p,p+1}_{1,1} (q)$ character is
related to the $M(p,p+1) +
\varphi_{2,1}$ model which is to be contrasted with the fermionic
representation (\ref{oneeight}) which is
related to $M(p,p+1) + \varphi_{1,3}$.\\

The object of this paper is to prove (\ref{oneten}) and its generalizations
which would, in
particular, provide a complete set of Virasoro characters in the limit $L
\rightarrow \infty$.
The rest of this article is organized as follows: In Section 2 we will
state a complete set of
polynomial identities proven in this paper.
In Section 3 we will discuss numerous properties of $q$ - trinomials and
derive recursion relations for
bosonic polynomials.\\
Section 4 is dedicated to deriving recurrences for fermionic polynomials by
means of a telescopic
expansion technique and proving identities stated in Section 2.\\
In Section 5 we take the limit $L \rightarrow \infty$ to obtain
Virasoro character identities.\\
We will conclude with some closing remarks.  Some technical details are
relegated to Appendices A and B.

\mysection{Identities}

Let us denote by $f^p_k (L; {\bf u}; {\bf A})$ the fundamental fermionic form
defined as
\bea
f^p_k (L; {\bf u}; {\bf A}) = \sum_{\scriptstyle {\bf m}\atop\scriptstyle
m_{p-2} \equiv k({\rm mod} \:
2)} q^{{{\bf m} \; {\tilde C}_{p-1} \; {\bf m} \over 4} - {{\bf A}
\cdot{\bf m} \over 2}}
\prod^{p-2}_{j=0} \[ \begin{array}{c} ({1 \over 2} {\tilde I}_{p-1} {\bf m}
+ L \; {\bf e}_0 + {1 \over 2}
{\bf u})_j \\ m_j \end{array} \]_q\nn
\label{twoone}
\ena
and for $a,b, i \in \zz_{\geq 0}$ let us set
\bea
{\bf u}^i_{a,b} = \theta (a > 1) {\bf e}_a + \theta (b>1) {\bf e}_b +
{\tilde p} (a+\delta_{a,1} + b+
\delta_{b,1} + \delta_{i,1} (p-1))({\bf e}_1 - {\bf e}_0)
\label{twotwo}
\ena
and
\bea
{\bf A}^i_{a,b} = \theta (b>1) {\bf e}_b + {\tilde p} (a+\delta_{a,1} + b +
\delta_{b,1} + \delta_{i,1} (p-1)) {\bf e}_1
\label{twothree}
\ena
with
\bea
\delta_{a,b} = \left\{ \begin{array}{c} 1 \; {\rm if} \; a=b \;  \\ 0 \;
{\rm if} \; a \not= b \; ,
\end{array}
\right.
\label{twofour}
\ena
\bea
{\tilde p} (a) = \left\{ \begin{array}{c} 1 \; {\rm if} \; a \equiv 1 \;
({\rm mod} \; 2 ) \\ 0 \;
{\rm if}
\; a \equiv 0 \; ({\rm mod} \; 2)  \; , \end{array} \right.
\label{twofive}
\ena
and
\bea
\theta \; ({\rm statement}) = \left\{ \begin{array}{l} 1 \; {\rm if \;
statement \; is \; true} \\ 0
\; {\rm if \; statement \; is \; false \; .}  \end{array} \right.
\label{twosix}
\ena
Next we define the objects of most interest to us
\bea
F^{p,i}_{a,b} (L) = \left\{ \begin{array}{l} f^p_i (L; {\bf u}^i_{a,b};
{\bf A}^i_{a,b} ) \;\;\:\: {\rm
if}
\; b \not= p-1 \\ f^p_{1-i} (L; {\bf u}^i_{a,b} ; {\bf A}_{a,b}) \; {\rm
if} \; b= p-1 \end{array}
\right.
\label{twoseven}
\ena
with $i=0,1$,
\bea
B^p_{a,b} (L,s)= \sum^\infty_{j=-\infty}  && \left\{ q^{pp'j^2 +j (p'a-ps)}
\( \begin{array}{c} L, q \\
2pj +a-b \end{array} \)^0 \right. \nn
&& \left. -q^{(pj +a)(p'j+s)} \( \begin{array}{c}  L, q \\ 2pj+a+b
\end{array} \)^0 \right\} \; ,
\label{twoeight}
\ena
and
\bea
{\tilde B}^p_{a,b} (L,s) = \sum^\infty_{j=-\infty} && \left\{ q^{pp'j^2 +j
(p'a-ps)} \( \begin{array}{c}
L, q \\ 2pj +a-b \end{array} \)^1 \right. \nn
&& \left. -q^{pp'j^2 +j(p'a+p(s-2)) +a(s-1)-b} \( \begin{array}{c}  L, q \\
2pj+a+b \end{array} \)^1
\right\}
\label{twonine}
\ena
with $p' = p+1$.
Finally, introducing functions $\varphi (a,b)$ and 
$ {\bf \tilde \varphi} (a,b)$
which are defined by
\bea
\varphi (a,b) = {(b-a)(b-a+1)   \over 4}                      
\label{twoten}
\ena
and
\bea
{\tilde \varphi} (a,b) = {(b-a)(b-a+3)\over 4} \; .
\label{twoeleven}
\ena
We are in position to state the main results of this paper.

\bigskip

\centerline{\underline{\large Even Identities}}
\bea
F^{p,0}_{a,1} (L) = \left\{ \begin{array}{cl} q^{{a(a-1) \over 4}}
B^p_{a,1} (L,1) & {\rm for} \; a+
\delta_{a,1} \equiv 0 \; ({\rm mod} \; 2) \\ q^{{(a-1)(a-2) \over 4}}
B^p_{a,1} (L,2) & {\rm for} \;
a+
\delta_{a,1} \equiv 1 \; ({\rm mod} \; 2) \end{array} \right.
\label{twotwelve}
\ena
\bea
F^{p,0}_{a,b} (L) = \left\{ \begin{array}{l} q^{\varphi (a,b)} B^p_{a,b}
(L, b+1) \; {\rm for} \; b+a +
\delta_{a,1} \equiv 0 \; ({\rm mod} \; 2) \\ q^{\tilde{\varphi} (a,b)}
{\tilde B}^p_{a,b} (L,b+2) \;\:
{\rm for} \; b+a + \delta_{a,1} \equiv 1 \; ({\rm mod} \; 2) \end{array}
\right.
\label{twothirteen}
\ena

\bigskip

\centerline{\underline{\large Odd Identities}}
\bea
F^{p,1}_{p-a,1} (L) = \left\{ \begin{array}{cl} q^{\varphi (p-a, p-1)}
B^p_{a,1} (L,1) &{\rm for} \; a
\equiv 1 \; ({\rm mod} \; 2) \\ q^{\tilde{\varphi}(p-a,p-1)+1-a}
B^p_{a,1} (L,2) & {\rm for} \; a \equiv
0 \; ({\rm mod} \; 2) \end{array} \right.
\label{twofourteen}
\ena
\bea
F^{p,1}_{p-a,b} (L) + \delta_{L,0} \delta_{a,b} = \left\{ \begin{array}{c}
q^{\varphi(p-a, p-1)+\varphi
(a,b) - {a-1 \over 4} a} B^p_{a,b} (L, b+1) \\ {\rm for} \; a+b \equiv 1 \;
({\rm mod} \; 2) \\
q^{{\tilde \varphi} (p-a,p-1)+ {\tilde \varphi}(a,b) - {a-1 \over 4} (a
+2)} {\tilde B}^p_{a,b}
(L, b+2) \\ {\rm for} \; a+b \equiv 0 \; ({\rm mod} \; 2) \end{array} \right.
\label{twofifteen}
\ena
with $F^{p,1}_{p-1, b} (L) \equiv F^{p,0}_{1,b} (L) \; {\rm for} \; 1 \leq
b \leq p-1$.
The identities (\ref{twotwelve}) - (\ref{twofifteen}) hold for $L=0, 1,2,
\ldots ; p= 4,5,6, \ldots ;
a=1,2, \ldots , p-2$ and $b=2,3, \ldots , p-1$.

Before we move on, a few comments are in order.  The polynomials $B_{a,b} (L,
b+1)$ and ${\tilde B}_{a,b} (L, b+2)$ (with
restrictions on $b$ - parameters which are not important for us here) were
introduced in \cite{bfifteen}-\cite{bsixteen} as configuration sums
for the regime $1^+$ of dilute $A_{p-1}$ models.  In addition, two more
polynomials $B^p_{a,b} (L,b)$ and
\bea
B^{'p}_{a,b} (L) &=& \sum^\infty_{j=-\infty} \left\{ q^{pp'j^2
+j(p'a-p(b+1))} \( \begin{array}{c} L,q \\
2pj +a-b \end{array} \)^1 \right. \nn
&-& \left. q^{pp'j^2 + j(p'a+p(b-1))+b(a-1)} \( \begin{array}{c} L,q \\
2pj+a+b \end{array} \)^1 \right\}
\label{twosixteen}
\ena
were employed in \cite{bsixteen} for the reason of completeness.
Here we point out
that the last two polynomials
are trivially related to the first two by
\bea
B^p_{p-a,p-b} (L, p-b+1) = B^p_{a,b} (L,b)
\label{twoseventeen}
\ena
\bea
{\tilde B}^p_{p-a, p-b} (L,p-b+2) = q^{b-a} B^{'p}_{a,b} (L) \; .
\label{twoeighteen}
\ena
In what follows we find it convenient to use simplified notations for
$F^{p,i}_{a,b} (L)$, namely
\bea
F^{p,i}_{a,b} (L) = \tilde{\sum}^{(i)} q^{\Phi^{p,i}_{a,b}
({\bf m})} \[ \begin{array}{c} {\bf n} + {\bf m} \\ {\bf n} \end{array}
\]^{p,L}_{a,b}
\label{twonineteen}
\ena
where for any $(p-1)$ - dimensional vectors ${\bf D}, {\bf E}$
\bea
\[ \begin{array}{c} {\bf n} + {\bf m} + {\bf D} \\ {\bf n} + {\bf E}
\end{array} \]^{p,L}_{a,b} =
\prod^{p-2}_{i=0} \[ \begin{array}{c} n_i + m_i +D_i \\ n_i + E_i
\end{array} \]_q
\label{twotwenty}
\ena
and $\tilde{\sum}^{(i)}$ stands for the sum over all ${\bf n}, {\bf m}$ subject
to constraints
\bea
{\bf n} + {\bf m} = {1 \over 2} {\tilde I}_{p-1} {\bf m} + L \; {\bf e}_0 +
{1 \over 2} {\bf u}^i_{a,b}
\label{twotwentyone}
\ena
and
\bea
m_{p-2} \equiv \left\{ \begin{array}{cl} i \; ({\rm mod} \; 2) & {\rm if}
\; b \not= p-1 \\ 1-i \; ({\rm
mod} \; 2) & {\rm if} \; b = p-1 \; , \end{array} \right.
\label{twotwentytwo}
\ena
and the quadratic form $\Phi^{p,i}_{a,b}({\bf m})$ is defined as
\bea
\Phi^{p,i}_{a,b} ({\bf m}) = {1 \over 4} {\bf m}  \tilde{C}_{p-1} {\bf m} -
{1 \over 2} {\bf
A}^i_{a,b} \cdot {\bf m} \; .
\label{twotwentythree}
\ena
We will also use a compact symbol
\bea
\{ {\bf n}, {\bf m} \}^{L,i}_{a,b}
\label{twotwentyfour}
\ena
to denote set of all integer solutions to the ${\bf n}, {\bf m}$ - system
(\ref{twotwentyone}) subject to the
additional constraint (\ref{twotwentytwo}).  (To keep notations tidy we
suppressed dependence on $p$ in
(\ref{twotwentyfour}).)  It immediately follows from (\ref{onesix}) and
(\ref{twonineteen}) that ${\bf
n}$ and ${\bf m}$ are vectors with non-negative integer components, i.e.
${\bf n}, {\bf m} \in
\zz^{p-1}_{\geq 0}$.\\
We now list a few easily verifiable consequences of (\ref{twotwentyone})
which will be frequently used.
\bea
L= n_0 + n_1 + 2 \sum^{p-3}_{j=1} j ({\bf n} - {1 \over 2} \; {\bf u}^i_{a,b}
)_{j+1} + (p-2) m_{p-2}\; ,
\label{twotwentyfive}
\ena
\bea
L= n_0 + n_1 + m_0 + m_1 \; ,
\label{twotwentysix}
\ena
\bea
L= n_0 + n_1 + m_2 + 2 \sum^{p-2}_{j=2} ({\bf n} - {1 \over 2} \; {\bf
u}^i_{a,b}
)_j +m_{p-2} \; ,
\label{twotwentyseven}
\ena
\bea
m_0 = \sum^{p-2}_{j=1} j ( {\bf n} - {1 \over 2} \; {\bf u}^i_{a,b})_j + {1
\over 2} (p-1) m_{p-2} \; ,
\label{twotwentyeight}
\ena
and
\bea
m_k = 2 \sum^{p-2-k}_{j=1} j ( {\bf n} - {1 \over 2} {\bf u}^i_{a,b}
)_{k+j} + (p-1-k) m_{p-2} \; {\rm
for} \; 2 \leq k \leq p-2 \; .
\label{twotwentynine}
\ena
In the spirit of Andrews proof of Rogers-Ramanujan identities
\cite{bseventeen} we will
prove
(\ref{twotwelve}-\ref{twofifteen}) by showing that fermionic and bosonic
polynomials there satisfy the
same recurrences and that our identities hold for $L =0$. \\
To generate recurrences for fermionic polynomials in
(\ref{twotwelve}-\ref{twofifteen}) we will employ
the technique of telescopic expansions which was introduced in
\cite{beleven} and further
developed in \cite{beighteen}-\cite{btwentythree}.  This
technique is based on the following elementary recursion relations for $q$
- binomial coefficients:
\bea
\[ \begin{array}{l} n+m \\ n \end{array} \]_q = \[ \begin{array}{l} n-1+m
\\ n-1 \end{array} \]_q + q^n
\[ \begin{array}{l} n+m-1 \\ n \end{array} \]_q
\label{twothirty}
\ena
and
\bea
\[ \begin{array}{l} n+m \\ n \end{array} \]_q = \[ \begin{array}{l} n+m-1
\\ n \end{array} \]_q + q^m
\[ \begin{array}{l} n-1+m \\ n-1 \end{array} \]_q
\label{twothirtyone}
\ena
Note, however, that (\ref{twothirty}), (\ref{twothirtyone}) fail for $n = m
= 0$.  In this case one
gets the contradiction $1=0$.  Fortunately, there is a way to modify the
definition (\ref{onesix}) in
such a way that ({\ref{twothirty}), (\ref{twothirtyone}) will hold without
exceptions.  The appropriate
definition which we borrow from \cite{btwentyone}  is
\bea
\[ \begin{array}{l} n+m\\ n \end{array} \]'_q = \left\{ \begin{array}{l}
{(q^{n+1})_m \over (q)_m} \;
{\rm for} \; m \in \zz_{\geq 0} , n \in \zz \\ 0 \quad - \; {\rm
otherwise} \end{array} \right.
\label{twothirtytwo}
\ena
It is easy to verify that modified $q$ - binomials (\ref{twothirtytwo}) can
be expressed in terms of
standard ones (\ref{onesix}) as
\bea
\[ \begin{array}{l} n+m \\ n \end{array} \]'_q = \left\{ \begin{array}{l}
\[ \begin{array}{l} n+m
\\ n \end{array} \]_q  \; {\rm if} \; n,m \geq 0 \; ,\\ (-)^m q^{{m+1+2n
\over 2} m} \[ \begin{array}{c}
-n-1
\\ m \end{array} \]_q \; {\rm if} \; n+m < 0 \; {\rm and} \;  m \geq 0 \\ 0
\quad - \; {\rm otherwise} .
\end{array}
\right.
\label{twothirtythree}
\ena
Clearly, modified $q$ - binomials vanish for negative values of $n$ as long
as $n+m \geq 0$.  At this
point one may wonder if it is possible to replace all standard binomials in
(\ref{twonineteen}) by
modified ones and then use (\ref{twothirty}) and (\ref{twothirtyone})
universally without hesitation.  Equation (\ref{twothirtythree}) indicates
that such a replacement
will be legitimate if one can
show that $n_i + m_i \geq 0$ for all $i=0,1,2, \ldots, p-2$.  From equation
(\ref{twothirtytwo}) we easily
infer that all $m$ - variables will remain non-negative, i.e.
\bea
m_i \geq 0 \; {\rm for} \; 0 \leq i \leq p-2 \; .
\label{twothirtyfour}
\ena
This being the case, it follows from (\ref{twotwentyone}) that $n_i + m_i
\geq 0$ for $i \not= 0$, which
in turn implies
\bea
n_i \geq 0 \; {\rm for} \; 1 \leq i \leq p-2
\label{twothirtyfive}
\ena
For $i=0$, the ${\bf n}, {\bf m}$ system (\ref{twotwentyone}) does not
guarantee
that $n_0+m_0 \geq 0$, and ,
in fact, the sum in (\ref{twonineteen}) will receive extra contributions from
negative values of $n_0$ if the
above mentioned replacement is made.  However, in Appendix A we will
demonstrate that these negative
contributions cancel each other out, unless $L=0, p-a=b \not= 1, p-1; \;
i=1$ in (\ref{twonineteen}).  In
the latter case one has an overall contribution equal to 1.  Thus, if
we replace all standard
binomials in (\ref{twonineteen}) by modified ones and denote the resulting
polynomials as $\bar{F}^{p,i}_{a,b} (L)$, then
\bea
\bar{F}^{p,0}_{a,b} (L) = F^{p,0}_{a,b} (L)
\label{twothirtysix}
\ena
and
\bea
\bar{F}^{p,1}_{a,b} (L) = F^{p,1}_{a,b} (L) + \delta_{L,0} \delta_{a,p-b}
\; \theta (1<b<p-1)
\label{twothirtyseven}
\ena
for $1 \leq a \leq p-2, \; 1 \leq b \leq p-1$.

\mysection{Bosonic Recurrences}

In the course of this paper we shall require various properties of $q$ -
trinomial coefficients such as
asymptotics, symmetry, recursion relations and some identities which become
tautologies when $q=1$.  Below
we shall first list necessary formulas and indicate where they first
appeared in the literature.  Next we
shall use the collected results to obtain two new $q$ - trinomial identities
and to derive recurrences for
bosonic polynomials introduced in Section 2.

\newpage

\centerline{\underline{\large Limiting Properties}}
\bea
\lim_{L \rightarrow \infty} \( \begin{array}{c} L, q \\ A \end{array} \)^0
= {1 \over (q)_\infty}
\label{threeone}
\ena
\bea
\lim_{L \rightarrow \infty} \( \begin{array}{c} L, q \\ A \end{array} \)^1
= {1+q^A \over
(q)_\infty}
\label{threetwo}
\ena
Equations (\ref{threeone}) and (\ref{threetwo}) are formulas (2.48) and
(2.49) of \cite{bone}, respectively.

\bigskip

\centerline{\underline{\large Symmetry Relations}}
\bea
\( \begin{array}{c} L, q \\ A \end{array} \)^0 = \( \begin{array}{c} L,q \\
-A \end{array} \)^0
\label{threethree}
\ena
\bea
\( \begin{array}{c} L, q \\ A \end{array} \)^1 = q^A \( \begin{array}{c}
L,q \\ -A \end{array} \)^1
\label{threefour}
\ena
The above two equations follow immediately from equation (\ref{twofifteen})
of \cite{bone}.

\bigskip

\centerline{\underline{\large Pascal triangle recurrences}}
\bea
\( \begin{array}{c} L,q \\ A \end{array} \)^0 = \( \begin{array}{c} L-1,q
\\ A+1 \end{array} \)^0 +
q^{L-1-A} \( \begin{array}{c} L-1,q \\A \end{array} \)^1 +q^{L-A} \(
\begin{array}{c} L-1,q \\ A-1
\end{array} \)^0
\label{threefive}
\ena
\bea
\( \begin{array}{c} L,q \\ A \end{array} \)^1 = \( \begin{array}{c} L-1,q
\\ A-1 \end{array} \)^0 +
q^{L-1} \( \begin{array}{c} L-1,q \\ A \end{array} \)^1 +q^A \(
\begin{array}{c} L-1,q \\ A+1
\end{array} \)^0
\label{threesix}
\ena
Depth one recurrences (\ref{threefive}) and (\ref{threesix}) were first
established in \cite{bone} where they
appeared as formulas (\ref{twotwentysix}) and (\ref{twotwentyfive}),
respectively.  In addition, depth
two recurrences
\bea
&& \( \begin{array}{c} L,q \\ A \end{array} \)^n = q^{L-A} \(
\begin{array}{c} L-1,q \\ A-1 \end{array}
\)^n + q^{L+A-n} \( \begin{array}{c} L-1,q \\ A+1 \end{array} \)^n + \nn [1em]
&& \quad + \( \begin{array}{c} L-1,q \\ A \end{array} \)^n + q^{L-1-n} \(
1-q^{L-1} \) \(
\begin{array}{c} L-2,q \\ A \end{array} \)^n
\label{threeseven}
\ena
with $n \in \zz$, have proven to be very useful.  For $n=0$, equation
(\ref{threeseven}) was first derived
in \cite{bfour} as equation (\ref{twothree}).  
Later, the proof given there was further
extended for all $n \in \zz$ in
\cite{beighteen} as equation (A.1).

\newpage

\centerline{\underline{\large Tautologies}}
\bea
\( \begin{array}{c} L,q \\ A \end{array} \)^1 = \( \begin{array}{c} L,q \\
A \end{array} \)^0 + q^A
(1-q^L) \( \begin{array}{c} L-1,q \\ A+1 \end{array} \)^0
\label{threeeight}
\ena
\bea
\( \begin{array}{c} L,q \\ A \end{array} \)^0 - q^{L-A} \( \begin{array}{c}
L,q \\ A \end{array} \)^1 =
\( \begin{array}{c} L,q \\ A+1 \end{array} \)^0 - q^L \( \begin{array}{c}
L,q \\ A+1 \end{array} \)^1
\label{threenine}
\ena
\bea
\( \begin{array}{c} L,q \\ A \end{array} \)^0 = q^L \( \begin{array}{c} L,q
\\ A \end{array} \)^1 &+&
(1-q^L) \( \begin{array}{c} L-1,q \\ A-1 \end{array} \)^0 \nn
&+& q^{L-1} (1-q^L) \( \begin{array}{c} L-1,q \\ A \end{array} \)^1
\label{threeten}
\ena
\bea
\( \begin{array}{c} L,q \\ A+1 \end{array} \)^1 - q^{1-A} \(
\begin{array}{c} L,q \\ A-1 \end{array} \)^1
= q^{A+1} \( \begin{array}{c} L-q \\ A+1 \end{array} \)^0 - q^{1-A} \(
\begin{array}{c} L,q \\ A-1
\end{array} \)^0
\label{threeeleven}
\ena
Identity (\ref{threeeight}) is equation (\ref{twotwentythree}) of
\cite{bone}  with
$A=B$.  Identity
({\ref{threenine}) is equation (\ref{twotwentyseven}) (corrected) of
\cite{bone}.
Identity (\ref{threeten}) can
be obtained from equation (\ref{twothirteen}) of \cite{btwentytwo}
 with $p=2, r=1$.
Finally, identity
(\ref{threeeleven}) follows from equation (A.8) of \cite{beighteen}
 with $n=1$ and $A$
replaced by $A+1$.

Each of identities (\ref{threefive})-(\ref{threeeleven}) gives rise to a
{\underline{\underline{twin}}} identity if
we first replace $A$ by $-A$ and then use symmetry relations
(\ref{threethree}), (\ref{threefour}).  In
particular, from (\ref{threefive}) we obtain
\bea
\( \begin{array}{c} L,q \\ A \end{array} \)^0 =  \( \begin{array}{c} L-1,q
\\ A-1 \end{array} \)^0 +
q^{L-1} \( \begin{array}{c} L-1,q \\ A \end{array} \)^1 + q^{L+A} \(
\begin{array}{c} L-1,q \\ A+1
\end{array} \)^0
\label{threetwelve}
\ena
Our object now is to prove two {\underline{new $q$ - trinomial coefficients
recurrences}}:
\bea
\( \begin{array}{c} L,q \\ A \end{array} \)^0 &=&  q^{L-A} \(
\begin{array}{c} L-1,q \\ A-1 \end{array}
\)^1 +  \( \begin{array}{c} L-1,q \\ A \end{array} \)^0 + q^{L-1} \(
\begin{array}{c} L-1,q \\ A+1
\end{array} \)^1\nn
&& \quad + q^{L-1} (q^{L-1} -1) \( \begin{array}{c} L-2,q \\ A \end{array} \)^0
\label{threethirteen}
\ena
\bea
\( \begin{array}{c} L,q \\ A \end{array} \)^1 &=&  \left\{ \(
\begin{array}{c} L-1,q \\ A-1 \end{array} \)^0 +
(q^{L-1}-1) \( \begin{array}{c} L-2,q \\ A-1 \end{array} \)^0 \right\} + \(
\begin{array}{c} L-1,q \\ A
\end{array} \)^1\nn
&+& q^A \left\{ \( \begin{array}{c} L-1,q \\ A+1 \end{array} \)^0 +
(q^{L-1} -1 ) \( \begin{array}{c}
L-2,q \\ A+1 \end{array} \)^0 \right\}\nn
&& \quad + q^{L-2} (q^{L-1} -1) \( \begin{array}{c} L-2,q \\ A \end{array}
\)^1 \; .
\label{threefourteen}
\ena
Both (\ref{threethirteen}) and (\ref{threefourteen}) are depth two
recurrences, but unlike relations
(\ref{threeseven}) they mix $q$ - trinomials with different superscripts.
The proof of (\ref{threethirteen})
follows immediately by combining (\ref{threeseven}) with $n=0$ and
(\ref{threeeight}) along with its
twin.  The proof of (\ref{threefourteen}) is just slightly more involved.
Here, we first combine
(\ref{threeten}) with the twin of (\ref{threenine}) to get
\bea
\( \begin{array}{c} L,q \\ A-1 \end{array} \)^0 &=&  q^{L+1-A} \(
\begin{array}{c} L,q \\ A-1
\end{array} \)^1 +  (1-q^L) \( \begin{array}{c} L-1,q \\ A-1 \end{array}
\)^0 \nn
&& \quad + q^{L-1} (1-q^L) \( \begin{array}{c} L-1,q \\ A \end{array} \)^1
\label{threefifteen}
\ena
and its twin
\bea
\( \begin{array}{c} L,q \\ A+1 \end{array} \)^0 &=&  q^L \(
\begin{array}{c} L,q \\ A+1
\end{array} \)^1 +  (1-q^L) \( \begin{array}{c} L-1,q \\ A+1 \end{array}
\)^0 \nn
&& \quad + q^{L-1-A} (1-q^L) \( \begin{array}{c} L-1,q \\ A \end{array}
\)^1 \; .
\label{threesixteen}
\ena
Next we use (\ref{threeseven}) with $n=1$ and (\ref{threefifteen}),
(\ref{threesixteen}) with $L$ replaced by $L-1$
to obtain (\ref{threefourteen}).

Now we are ready to establish the following recursion relations for bosonic
polynomials introduced in
Section 2:
\bea
B^p_{a,1} (L,1) = B^p_{a,1} (L-1, 1) + q^{L+1-a} B^p_{a,2} (L-1,3)
\label{threeseventeen}
\ena
\bea
B^p_{a,1} (L,2) = B^p_{a,1} (L-1, 2) + q^{L+1-a} \tilde{B}^p_{a,2} (L-1,4)
\label{threeeighteen}
\ena
\bea
B^p_{a,2} (L,3) = B^p_{a,2} (L-1, 3) &+& q^{L+a-2} B^p_{a,1} (L-1,1)\nn
&+& q^{L-a+2} \tilde{B}^p_{a,3} (L-1,5)
\label{threenineteen}
\ena
\bea
\tilde{B}^p_{a,2} (L,4) = q^{L-1} \tilde{B}_{a,2} (L-1, 4) &+& q^{a-2}
B^p_{a,1} (L-1,2)\nn
&+& B^p_{a,3} (L-1,4)
\label{threetwenty}
\ena
\bea
B^p_{a,b} (L,b+1) &=& B^p_{a,b} (L-1, b+1) + q^{L-1} (q^{L-1} -1) B^p_{a,b}
(L-2,b+1)\nn
&+& q^{L-a+b} \tilde{B}^p_{a,b+1} (L-1,b+3) +q^{L-1} \tilde{B}^p_{a,b-1}
(L-1, b+1)
\label{threetwentyone}
\ena
\bea
\tilde{B}^p_{a,b} (L,b+2) &=& B^p_{a,b+1} (L-1, b+2) + (q^{L-1}-1)
B^p_{a,b+1} (L-2,b+2)\nn
&+& q^{a-b} (B^p_{a,b-1} (L-1,b) + (q^{L-1}-1) B^p_{a,b-1} (L-2, b))\nn
&+& \tilde{B}^p_{a,b} (L-1, b+2) + q^{L-2} (q^{L-1}-1)  \tilde{B}^p_{a,b}
(L-2, b+2)
\label{threetwentytwo}
\ena
Because the factor $q^{L-1}-1$ in (\ref{threetwentyone}) and
(\ref{threetwentytwo}) vanishes for $L=1$ one needs to specify initial
conditions only for $L=0$ to determine these polynomials uniquely.
The validity of (\ref{threeseventeen})-(\ref{threenineteen}) is readily
seen from (\ref{threefive}) and
(\ref{threetwelve}).  Relation (\ref{threetwenty}) follows by use of
(\ref{threesix}) and its twin
identity.

Finally (\ref{threethirteen}) along with its twin and (\ref{threefourteen})
along with its twin imply
(\ref{threetwentyone}) and (\ref{threetwentytwo}), respectively.  We close
this section by noting some
additional properties of bosonic polynomials.

\bigskip

\centerline{\underline{\large Reflection properties}}
\bea
B^p_{a,p-1} (L,p) = B^p_{p-a,1} (L,1)
\label{threetwentythree}
\ena
\bea
\tilde{B}^p_{a,p-1} (L,p+1) = q^{a-p+1} B^p_{p-a,1} (L,2)
\label{threetwentyfour}
\ena

\bigskip

\centerline{\underline{\large Limiting properties}}
\bea
\lim_{L \rightarrow \infty} B^p_{a,b} (L,b+1) = \chi^{p,p+1}_{a,b+1} (q)
\label{threetwentyfive}
\ena
\bea
\lim_{L \rightarrow \infty} \tilde{B}^p_{a,b} (L,b+2) =
\chi^{p,p+1}_{a,b+2} (q) + q^{a-b} \chi^{p,p+1}_{a,b}
(q)
\label{threetwentysix}
\ena
with $\chi^{p,p+1}_{a,b} (q)$ defined by (\ref{oneseven}).

\bigskip

\centerline{\underline{\large Closing properties}}
\bea
B^p_{a,p} (L, p+1) = \tilde{B}^p_{a,p} (L, p+2) =0
\label{threetwentyseven}
\ena
Equations (\ref{threetwentythree}) and (\ref{threetwentyseven}) can be
easily verified by inspection.
Equations (\ref{threetwentyfive})-(\ref{threetwentysix}) are simply
consequences of (\ref{threeone}) and
(\ref{threetwo}).  Identity (\ref{threetwentyfour}) follows from
(\ref{threeeleven}).  Tautology
(\ref{threeeleven}) can be also used to derive a further interesting identity
\bea
B^p_{a,1} (L, 2) = \tilde{B}^p_{a,1} (L, 2) \; .
\label{threetwentyeight}
\ena

\mysection{Fermionic Recursion Relations}

We now turn our attention to the proof of the following fermionic recurrences:
\bea
F^{p,0}_{a,1} (L) = F^{p,0}_{a,1} (L-1) + q^{L -\frac12} F^{p,0}_{a,2}
(L-1) \; {\rm for} \; a+ \delta_{a,1}
\equiv 0 \; ({\rm mod} \; 2)
\label{fourone}
\ena
\bea
F^{p,0}_{a,1} (L) = F^{p,0}_{a,1} (L-1) + q^{L -1} F^{p,0}_{a,2} (L-1) \;
{\rm for} \; a+ \delta_{a,1}
\equiv 1 \; ({\rm mod} \; 2)
\label{fourtwo}
\ena
\bea
F^{p,0}_{a,2} (L) &=& q^{L -1} F^{p,0}_{a,3} (L-1) +q^{L -\frac12}
F^{p,0}_{a,1} (L-1) + F^{p,0}_{a,2} (L-1) \nn
& & \qquad {\rm for} \; a + \delta_{a,1} \equiv 0 \; ({\rm mod} \; 2)
\label{fourthree}
\ena
\bea
F^{p,0}_{a,2} (L) &=& {1 \over \sqrt{q}} F^{p,0}_{a,3} (L-1) +
F^{p,0}_{a,1} (L-1) + q^{L -1}F^{p,0}_{a,2} (L-1) \nn
& & \qquad {\rm for} \; a + \delta_{a,1} \equiv 1 \; ({\rm mod} \; 2)
\label{fourfour}
\ena
\bea
F^{p,0}_{a,b} (L) &=& q^{L-1} F^{p,0}_{a,b+1} (L-1) + q^{L-\frac12}
F^{p,0}_{a,b-1} (L-1) + F^{p,0}_{a,b} (L-1) \nn
&+& q^{L-1} (q^{L-1} -1) F^{p,0}_{a,b} (L-2) \nn
& & \qquad {\rm for} \; b+a+ \delta_{a,1} \equiv 0 \; ({\rm mod} \; 2), 2<b
\leq p-1
\label{fourfive}
\ena
\bea
F^{p,0}_{a,b} (L) &=& {1 \over \sqrt{q}} \{ F^{p,0}_{a,b+1} (L-1) +
(q^{L-1}-1) F^{p,0}_{a,b+1} (L-2)
\} \nn
&+& \left\{ F^{p,0}_{a,b-1} (L-1) + (q^{L-1} -1) F^{p,0}_{a,b-1}  (L-2)
\right\} \nn
&+& \left\{ F^{p,0}_{a,b} (L-1) + q^{L-2} (q^{L-1} -1) F^{p,0}_{a,b}  (L-2)
\right\} \nn
& & \qquad {\rm for} \; b+a+ \delta_{a,1} \equiv 1 \; ({\rm mod} \; 2), 2<b
\leq p-1
\label{foursix}
\ena
with $F^{p,0}_{a,p} (L) \equiv 0$.
In contrast with many identities relating $q$ - trinomial coefficients used
in the proof of the bosonic recurrences
(3.17-22) the only identities we shall require here are the elementary
recursion relations for $q$ - binomials
(2.30-2.31).  Recurrences (4.1-5) will be proven in this section.  The
proof of equation (4.6) will be relegated to
Appendix B.

\bigskip

\centerline{\underline{\large 4.1 Proof of (4.1) and (4.2)}}
\bigskip
Using equation (2.30) we can expand $F^{p,0}_{a,1} (L)$ as
\bea
F^{p,0}_{a,1} (L) &=& \widetilde{\sum}^{(0)} q^{\Phi^{p,0}_{a,1} ({\bf m})}
\[ \begin{array}{l} {\bf n} + {\bf m} -
{\bf e}_0 \\ {\bf n} - {\bf e}_0 \end{array} \]^{p,L}_{a,1} +\nn
&+& \widetilde{\sum}^{(0)} q^{\Phi^{p,0}_{a,1} ({\bf m}) + n_0} \[
\begin{array}{l} {\bf n} + {\bf m} - {\bf
e}_0 \\ {\bf n}\end{array} \]^{p,L}_{a,1} \; .
\label{fourseven}
\ena
Now equations (2.25) and (2.28-29) imply that
\bea
\{ {\bf n}, {\bf m} \}^{L,0}_{a,1} - \{ {\bf e}_0, 0 \} = \{ {\bf n}, {\bf
m} \}^{L-1, 0}_{a,1}
\label{foureight}
\ena
and
\bea
\{ {\bf n}, {\bf m} \}^{L,0}_{a,1} - \{ 0, {\bf e}_0 \} = \{ {\bf n}, {\bf
m} \}^{L-1, 0}_{a,2}
\label{fournine}
\ena
The first equation of (2.21)
\bea
n_0 + m_0 = L - {m_2 \over 2} - {1 \over 2} \; \tilde{p} (a+ \delta_{a,1}
+b+ \delta_{b,1})
\label{fourten}
\ena
with $b=1$ gives
\bea
\Phi^{p,0}_{a,1} ({\bf m}) + n_0 = L - {1 \over 2} - {1 \over 2} \tilde{p}
\; (a+ \delta_{a,1}) + \Phi^{p,0}_{a,2}
({\bf m} - {\bf e}_0) \; .
\label{foureleven}
\ena
Hence, we can identify the first and second terms in the right-hand-side of
(4.7) as \newline $F^{p,0}_{a,1} (L-1)$
and
$q^{L- \frac12 - \frac12 \tilde{p} (a + \delta_{a,1})} F^{p,0}_{a,2}
(L-1)$, respectively.  Thus (4.1) and (4.2) are
established.

\bigskip

\centerline{\underline{\large 4.2 Proof of (4.3) and (4.4)}}
\bigskip
To prove (4.3), we start by expanding $F^{p,0}_{a,2} (L)$ with
$a+ \delta_{a,1} \equiv 0 ({\rm mod} \; 2)$
in telescopic fashion as
\bea
F^{p,0}_{a,2} (L) &=& \widetilde{\sum}^{(0)} q^{\Phi^{p,0}_{a,2} ({\bf m})
+n_0} \[ \begin{array}{l} {\bf n} + {\bf
m} - {\bf e}_0 - {\bf e}_2 \\ {\bf n} \end{array} \]^{p,L}_{a,2} \nn
&+& \widetilde{\sum}^{(0)} q^{\Phi^{p,0}_{a,2} ({\bf m}) + n_0+m_2} \[
\begin{array}{l} {\bf n} + {\bf m} - {\bf
e}_0  - {\bf e}_2 \\ {\bf n} - {\bf e}_2 \end{array} \]^{p,L}_{a,2} \nn
&+& \widetilde{\sum}^{(0)} q^{\Phi^{p,0}_{a,2} ({\bf m})} \[
\begin{array}{l} {\bf n} + {\bf m} - {\bf
e}_0  \\ {\bf n}-  {\bf e}_0 \end{array} \]^{p,L}_{a,2}
\label{fourtwelve}
\ena
where (2.31) was used to combine first and second terms in the
right-hand-side of (4.12) and (2.30) was used in the
last step.

Next with the help of (2.25) and (2.28-29) one shows that
\bea
\{ {\bf n}, {\bf m} \}^{L,0}_{a,2} - \{ 0, {\bf e}_0 + {\bf e}_2  \} = \{
{\bf n}, {\bf m} \}^{L-1, 0}_{a,3} \; ,
\label{fourthirteen}
\ena
\bea
\{ {\bf n}, {\bf m} \}^{L,0}_{a,2} - \{ {\bf e}_2, {\bf e}_0  \} = \{ {\bf
n}, {\bf m} \}^{L-1, 0}_{a,3} \; ,
\label{fourfourteen}
\ena
and
\bea
\{ {\bf n}, {\bf m} \}^{L,0}_{a,2} - \{  {\bf e}_0, 0 \}  = \{ {\bf n},
{\bf m} \}^{L-1, 0}_{a,2}\; .
\label{fourfiftteen}
\ena
Using equation (4.10) with $b=2$ and $a+ \delta_{a,1} \equiv 0 ({\rm mod}
\; 2)$ one readily checks that
\bea
\Phi^{p,0}_{a,2} ({\bf m}) + n_0 = L-1 +  \Phi^{p,0}_{a,3} ({\bf m} - {\bf
e}_0 - {\bf e}_2)
\label{foursixteen}
\ena
and
\bea
\Phi^{p,0}_{a,2} ({\bf m}) + n_0 + m_2 = L- {1 \over 2} + \Phi^{p,0}_{a,1}
({\bf m} - {\bf e}_0) \; .
\label{fourseventeen}
\ena
Equations (4.13-17) enable one to recognize first, second, and third terms
in the right-hand-side of (4.12) as
$q^{L-1}  F^{p,0}_{a,3} (L-1), q^{L-\frac12} F^{p,0}_{a,1} (L-1)$, and
$F^{p,0}_{a,2} (L-1)$, respectively.
Thus, (4.3) is established.

The proof of (4.4) is however
considerably more involved then the simple demonstration of (4.3).  Again,
we start by
expanding $F^{p,0}_{a,2} (L)$ with $a+ \delta_{a,1} \equiv 1 ({\rm mod} \;
2)$ in telescopic fashion as
\bea
F^{p,0}_{a,2} (L) &=& \widetilde{\sum}^{(0)} q^{\Phi^{p,0}_{a,b} ({\bf m})}
\[ \begin{array}{l} {\bf n} + {\bf
m} - {\bf e}_1 - {\bf e}_2 \\ {\bf n} \end{array} \]^{p,L}_{a,b} +\nn
&+& \widetilde{\sum}^{(0)} q^{\Phi^{p,0}_{a,b} ({\bf m}) + m_2} \[
\begin{array}{l} {\bf n} + {\bf m} - {\bf
e}_1 - {\bf e}_2 \\ {\bf n} - {\bf e}_2 \end{array} \]^{p,L}_{a,b} \nn
&+& \widetilde{\sum}^{(0)} q^{\Phi^{p,0}_{a,2} ({\bf m}) + m_1} \[
\begin{array}{l} {\bf n} + {\bf m} - {\bf
e}_1  \\ {\bf n}-  {\bf e}_1 \end{array} \]^{p,L}_{a,b}
\label{foureighteen}
\ena
where (2.31) was used twice to recombine all terms in the right-hand-side.
In the same manner as before, it is
possible to show that
\bea
{1 \over \sqrt{q}} F^{p,0}_{a,3} (L-1) = \widetilde{\sum}^{(0)}
q^{\Phi^{p,0}_{a,2} ({\bf m})} \[
\begin{array}{l} {\bf n} + {\bf m} - {\bf e}_1 - {\bf e}_2 \\ {\bf n}
\end{array} \]^{p,L}_{a,b} \; ,
\label{fournineteen}
\ena
\bea
F^{p,0}_{a,1} (L-1) = \widetilde{\sum}^{(0)} q^{\Phi^{p,0}_{a,2} ({\bf m})
-m_0 + {m_2 +1 \over 2}} \[
\begin{array}{l} {\bf n} + {\bf m} - {\bf e}_0 - {\bf e}_2 \\ {\bf n} -
{\bf e}_2 \end{array} \]^{p,L}_{a,2} \; ,
\label{fourtwenty}
\ena
and
\bea
F^{p,0}_{a,2} (L-1) = \widetilde{\sum}^{(0)} q^{\Phi^{p,0}_{a,2} ({\bf m})} \[
\begin{array}{l} {\bf n} + {\bf m} - {\bf e}_0  \\ {\bf n} - {\bf e}_0
\end{array} \]^{p,L}_{a,2} \; .
\label{fourtwentyone}
\ena

Employing (2.31) one can expand (4.20) further as
\bea
F^{p,0}_{a,1} (L-1) &=& \widetilde{\sum}^{(0)} q^{\Phi^{p,0}_{a,2} ({\bf
m})-m_0+m_1+ {m_2 +1 \over 2}}  \[
\begin{array}{l} {\bf n} + {\bf m} - {\bf E}_{0,2}  \\ {\bf n} - {\bf
E}_{1,2}\end{array} \]^{p,L}_{a,2} \nn
&+& \widetilde{\sum}^{(0)} q^{\Phi^{p,0}_{a,2} ({\bf m}) - m_0 + {m_2 +1
\over 2}}  \[ \begin{array}{l} {\bf n} +
{\bf m} - {\bf E}_{0,2} \\ {\bf n} - {\bf e}_2 \end{array} \]^{p,L}_{a,2}
\label{fourtwentytwo}
\ena
where for $b \geq a \;\; (p-1)$ - dimensional vector ${\bf E}_{a,b}$ is
defined as
\bea
{\bf E}_{a,b} = \sum^b_{i=a} {\bf e}_i \; .
\label{fourtwentythree}
\ena
Next making use of (2.30) we expand the second term in the right-hand-side
of (4.18) as
\bea
\widetilde{\sum}^{(0)} q^{\Phi^{p,0}_{a,2} ({\bf m}) +m_2}  \[
\begin{array}{l} {\bf n} + {\bf m} - {\bf E}_{1,2}  \\ {\bf n} - {\bf e}_2
\end{array} \]^{p,L}_{a,2} &=&
\widetilde{\sum}^{(0)} q^{\Phi^{p,0}_{a,2} ({\bf m}) + m_2} \[
\begin{array}{l} {\bf n} + {\bf m} - {\bf E}_{0,2} \\
{\bf n} - {\bf e}_0 - {\bf e}_2 \end{array} \]^{p,L}_{a,2} + \nn
&+& \widetilde{\sum}^{(0)} q^{\Phi^{p,0}_{a,2} ({\bf m}) + m_2 + n_0} \[
\begin{array}{l} {\bf n} +
{\bf m} - {\bf E}_{0,2} \\ {\bf n} - {\bf e}_2 \end{array} \]^{p,L}_{a,2}
\label{fourtwentyfour}
\ena
and then perform the change of summation variables
\bea
{\bf n} \rightarrow {\bf n} + {\bf e}_0 - {\bf e}_1 + {\bf e}_2 \; ; {\bf
m} \rightarrow {\bf m} - {\bf e}_0 + {\bf
e}_1
\label{fourtwentyfive}
\ena
in the first sum appearing in the right-hand-side of (4.24) to obtain
\bea
\widetilde{\sum}^{(0)} q^{\Phi^{p,0}_{a,2} ({\bf m}) +m_2}  \[
\begin{array}{l} {\bf n} + {\bf m} - {\bf E}_{0,2}  \\ {\bf n} - {\bf e}_0
- {\bf e}_2 \end{array} \]^{p,L}_{a,2} =
\widetilde{\sum}^{(0)} q^{\Phi^{p,0}_{a,b} ({\bf m})- m_0 +m_1+
\frac{m_2 +1}{2}} \[ \begin{array}{l} {\bf n} + {\bf
m} - {\bf E}_{0,2} \\ {\bf n} - {\bf E}_{1,2} \end{array} \]^{p,L}_{a,2} \; .
\nn &&
\label{fourtwentysix}
\ena
Combining (4.22), (4.24), (4.26) and making use of (4.10) with $b=2, a+
\delta_{a,1} \equiv 1 ({\rm mod} \; 2)$, we
have
\bea
&& \widetilde{\sum}^{(0)} q^{\Phi^{p,0}_{a,2} ({\bf m}) +m_2}  \[
\begin{array}{l} {\bf n} + {\bf m} - {\bf E}_{1,2}  \\ {\bf n} - {\bf e}_2
\end{array} \]^{p,L}_{a,2} =\nn
&& = F_1 (L-1) + (q^{L-1} -1) \widetilde{\sum}^{(0)} q^{\Phi^{p,0}_{a,2}
({\bf m})- m_0 + {m_2 +1 \over 2}} \[
\begin{array}{l} {\bf n} + {\bf m} - {\bf E}_{0,2} \\ {\bf n} - {\bf e}_2
\end{array} \]^{p,L}_{a,2} \; .
\label{fourtwentyseven}
\ena
Finally, let us perform the change of summation variables
\bea
{\bf n} \rightarrow {\bf n} - {\bf e}_0 - {\bf e}_1 + {\bf e}_2  \; ;
{\bf m} \rightarrow {\bf m} + {\bf e}_0 + {\bf e}_1
\label{fourtwentyeight}
\ena
in the second sum in (4.27) to derive for the second term in the
right-hand-side of (4.18)
\bea
&& \widetilde{\sum}^{(0)} q^{\Phi^{p,0}_{a,2} ({\bf m}) +m_2}  \[
\begin{array}{l} {\bf n} + {\bf m} - {\bf E}_{1,2}  \\ {\bf n} - {\bf e}_2
\end{array} \]^{p,L}_{a,2} =\nn
&&  F_1 (L-1) + (q^{L-1} -1) \widetilde{\sum}^{(0)} q^{\Phi^{p,0}_{a,2}
({\bf m}) +m_1}  \[
\begin{array}{l} {\bf n} + {\bf m} - {\bf E}_{0,1}  \\ {\bf n} - {\bf
E}_{0,1} \end{array} \]^{p,L}_{a,2} \; .
\label{fourtwentynine}
\ena
We now turn to the last term in the right-hand-side of (4.18).  We start by
expanding it with the aid of (2.30) as
\bea
\widetilde{\sum}^{(0)} q^{\Phi^{p,0}_{a,2} ({\bf m}) +m_1}  \[
\begin{array}{l} {\bf n} + {\bf m} - {\bf e}_1  \\ {\bf n} - {\bf e}_1
\end{array} \]^{p,L}_{a,2} &=&
\widetilde{\sum}^{(0)} q^{\Phi^{p,0}_{a,2} ({\bf m}) + m_1 + n_0} \[
\begin{array}{l} {\bf n} + {\bf m} - {\bf E}_{0,1}
\\ {\bf n} - {\bf e}_1  \end{array} \]^{p,L}_{a,2}  \nn
&+& \widetilde{\sum}^{(0)} q^{\Phi^{p,0}_{a,2} ({\bf m}) + m_1} \[
\begin{array}{l} {\bf n} +
{\bf m} - {\bf E}_{0,1} \\ {\bf n} - {\bf E}_{0,1} \end{array}
\]^{p,L}_{a,2} \; .
\label{fourthirty}
\ena
Next we perform the charge of summation variables
\bea
{\bf n} \rightarrow {\bf n} - {\bf e}_0 + {\bf e}_1 \; ; {\bf m}
\rightarrow {\bf m} + {\bf e}_0 - {\bf
e}_1
\label{fourthirtyone}
\ena
in the first term in the right-hand-side of (4.30) to get
\bea
\widetilde{\sum}^{(0)} q^{\Phi^{p,0}_{a,2} ({\bf m}) + m_1}  \[
\begin{array}{l} {\bf n} + {\bf m} - {\bf e}_1  \\ {\bf n} - {\bf e}_1
\end{array} \]^{p,L}_{a,2}  &=&
\widetilde{\sum}^{(0)} q^{\Phi^{p,0}_{a,2} ({\rm m}) + L-1} \[
\begin{array}{l} {\bf n} + {\bf m} - {\bf E}_{0,1} \\
{\bf n} - {\bf e}_0  \end{array} \]^{p,L}_{a,2}  \nn
&+& \widetilde{\sum}^{(0)} q^{\Phi^{p,0}_{a,2} ({\rm m})+ m_1} \[
\begin{array}{l} {\bf n} +
{\bf m} - {\bf E}_{0,1} \\ {\bf n} - {\bf E}_{0,1} \end{array} \]^{p,L}_{a,2}
\label{fourthirtytwo}
\ena
where we also used (4.10) with $b=2, a + \delta_{a,1} \equiv 1 ({\rm mod}
\; 2)$.
If we now add and subtract from the right-hand-side of (4.32) the
following sum
\begin{eqnarray*}
\widetilde{\sum}^{(0)} q^{\Phi^{p,0}_{a,2} ({\rm m}) +m_1 + (L-1)} \[
\begin{array}{l} {\bf n} + {\bf m} - {\bf
E}_{0,1} \\ {\bf n} - {\bf E}_{0,1}  \end{array} \]^{p,L}_{a,2}
\end{eqnarray*}
and then use (2.31) to recombine two terms, the result is
\bea
&& \widetilde{\sum}^{(0)} q^{\Phi^{p,0}_{a,2} ({\bf m}) + m_1}  \[
\begin{array}{l} {\bf n} + {\bf m} - {\bf
e}_1 \\ {\bf n} - {\bf e}_1  \end{array} \]^{p,L}_{a,2} = \nn
&& = \widetilde{\sum}^{(0)} q^{\Phi^{p,0}_{a,2} ({\bf m}) + L-1} \[
\begin{array}{l} {\bf n} + {\bf m} - {\bf
e}_0 \\ {\bf n} - {\bf e}_0  \end{array} \]^{p,L}_{a,2} + \nn
&& + (1 - q^{L-1})  \widetilde{\sum}^{(0)} q^{\Phi^{p,0}_{a,2} ({\bf m}) +
m_1} \[ \begin{array}{l} {\bf n} + {\bf m}
- {\bf E}_{0,1} \\ {\bf n} - {\bf E}_{0,1}  \end{array} \]^{p,L}_{a,2} \; .
\label{fourthirtythree}
\ena
Recalling (4.21), we can rewrite (4.33) as
\bea
&& \widetilde{\sum}^{(0)} q^{\Phi^{p,0}_{a,2} ({\bf m}) + m_1}  \[
\begin{array}{l} {\bf n} + {\bf m} - {\bf
e}_1 \\ {\bf n} - {\bf e}_1  \end{array} \]^{p,L}_{a,2} = \nn
&& = q^{L-1} F^{p,0}_{a,2} (L-1) + (1-q^{L-1}) \widetilde{\sum}^{(0)}
q^{\Phi^{p,0}_{a,2} ({\bf m}) + m_1} \[
\begin{array}{l} {\bf n} + {\bf m} - {\bf E}_{0,1} \\ {\bf n} - {\bf
E}_{0,1}  \end{array} \]^{p,L}_{a,2} \; .
\label{fourthirtyfour}
\ena
Finally, combining (4.18), (4.19), (4.29) and (4.34) we arrive at the
desired result (4.4).  Thus (4.4) is
established.

\bigskip

\centerline{\underline{\large 4.3 Proof of (4.5) }}
\bigskip
Once again, we start by using (2.30), (2.31) to expand $F^{p,0}_{a,0} (L)$
with $a+ \delta_{a,1} +b+ \delta_{b,1}
\equiv 0 ({\rm mod} \; 2), 2<b \leq p-1$ in telescopic fashion as
\bea
F^{p,0}_{a,b} (L) &=& \theta (b \not= p-1) \widetilde{\sum}^{(0)}
q^{\Phi^{p,0}_{a,b} ({\bf m}) +n_0} \[
\begin{array}{l} {\bf n} + {\bf m} - {\bf e}_0 - {\bf E}_{2,b} \\ {\bf n}
\end{array} \]^{p,L}_{a,b} +\nn
&+& \sum^b_{l=3} \widetilde{\sum}^{(0)} q^{\Phi^{p,0}_{a,b} ({\bf m}) + n_0
+m_l} \[ \begin{array}{l} {\bf n} +
{\bf m} - {\bf e}_0 - {\bf E}_{2,l} \\ {\bf n} - {\bf e}_L \end{array}
\]^{p,L}_{a,b} \nn
&+& \widetilde{\sum}^{(0)} q^{\Phi^{p,0}_{a,b} ({\bf m}) + n_0 + m_2} \[
\begin{array}{l} {\bf n} + {\bf m} - {\bf
e}_0  - {\bf e}_2 \\ {\bf n}-  {\bf e}_2 \end{array} \]^{p,L}_{a,b} \nn
&+& \widetilde{\sum}^{(0)} q^{\Phi^{p,0}_{a,b} ({\bf m}) } \[
\begin{array}{l} {\bf n} + {\bf m} - {\bf
e}_0   \\ {\bf n}-  {\bf e}_0 \end{array} \]^{p,L}_{a,b} \; .
\label{fourthirtyfive}
\ena
Employing (2.25), (2.28), (2.29) along with (4.10) it is
straightforward (through tedious) to verify the following:
\bea
q^{L-1} F^{p,0}_{a,b+1} (L-1) = \theta (b \not= p-1) \widetilde{\sum}^{(0)}
q^{\Phi^{p,0}_{a,b} ({\bf m}) +n_0} \[
\begin{array}{l} {\bf n} + {\bf m} - {\bf e}_0 - {\bf E}_{2,b} \\ {\bf n}
\end{array} \]^{p,L}_{a,b}
\label{fourthirtysix}
\ena
\bea
F^{p,0}_{a,b} (L-1) = \widetilde{\sum}^{(0)} q^{\Phi^{p,0}_{a,b} ({\bf m})}
\[ \begin{array}{l}
{\bf n} + {\bf m} - {\bf e}_0 \\ {\bf n} - {\bf e}_0 \end{array}
\]^{p,L}_{a,b}
\label{fourthirtyseven}
\ena
\bea
q^{L-\frac12} F^{p,0}_{a,b-1} (L-1) = \widetilde{\sum}^{(0)}
q^{\Phi^{p,0}_{a,b} ({\bf m}) + n_0 +1 - m_{b-1} +m_b}
\[ \begin{array}{l} {\bf n} + {\bf m} - {\bf e}_0 - {\bf E}_{2,b-1} + {\bf
e}_{b-1} - {\bf e}_b
\\ {\bf n} +  {\bf e}_{b-1} - {\bf e}_b \end{array} \]^{p,L}_{a,b} \nn &&
\label{fourthirtyeight}
\ena
\bea
q^{L-1} F^{p,0}_{a,b} (L-2) = \widetilde{\sum}^{(0)} q^{\phi^{p,0}_{a,b}
({\bf m}) + n_0 + 1-m_{b-1} + m_b} \[
\begin{array}{l} {\bf n} + {\bf m} - {\bf E}_{0,1} - 2 {\bf E}_{2,b-1} +
{\bf e}_{b-1} - {\bf e}_b \\ {\bf n} +  {\bf
e}_{b-1} - {\bf e}_b \end{array} \]^{p,L}_{a,b} \nn &&
\label{fourthirtynine}
\ena
\bea
q^{2 (L-1)} F^{p,0}_{a,b} (L-2) = \widetilde{\sum}^{(0)}
q^{\Phi^{p,0}_{a,b} ({\bf m}) + 2n_0 + m_2} \[
\begin{array}{l} {\bf n} + {\bf m} - 2 {\bf e}_0 - {\bf e}_2 \\ {\bf n}-
{\bf e}_2 \end{array} \]^{p,L}_{a,b}
\label{fourforty}
\ena
with $m_{p-1} \equiv 0$ and ${\bf e}_{p-1} \equiv 0$.

To proceed further, it's expedient to expand (4.38) with the aid of (2.31)
in telescopic fashion as follows
\bea
&& q^{L-\frac12} F^{p,0}_{a,b-1} (L-1) = \widetilde{\sum}^{(0)}
q^{\Phi^{p,0}_{a,b} ({\bf m}) + n_0 + 1-m_{b-1} +
m_b} \[ \begin{array}{l} {\bf n} + {\bf m} - {\bf e}_0  - 2 {\bf E}_{2,b-1}
+ {\bf e}_{b-1} - {\bf e}_b \\ {\bf n} +
{\bf e}_{b-1} - {\bf e}_b \end{array} \]^{p,L}_{a,b}  \nn
&& + \sum^{b-1}_{l=2}  \widetilde{\sum}^{(0)} q^{\Phi^{p,0}_{a,b} ({\bf m})
+ n_0 + 1-m_{b-1} +
m_b + (m_l -1)} \[ \begin{array}{l} {\bf n} + {\bf m} - {\bf e}_0  - 2 {\bf
E}_{2,b-1} + {\bf e}_{b-1} - {\bf e}_b - {\bf
E}_{l,b-1} \\ {\bf n} +  {\bf e}_{b-1} - {\bf e}_b - {\bf e}_l \end{array}
\]^{p,L}_{a,b} \; . \nn &&
\label{fourfortyone}
\ena
Using (2.31) again, we split the first sum in (4.41) into two pieces as
\bea
&& \widetilde{\sum}^{(0)} q^{\Phi^{p,0}_{a,b} ({\bf m}) + n_0 + 1-m_{b-1} +
m_b} \[ \begin{array}{l} {\bf n} + {\bf m} - {\bf e}_0  - 2 {\bf E}_{2,b-1}
+ {\bf e}_{b-1} - {\bf e}_b \\
{\bf n} +  {\bf e}_{b-1} - {\bf e}_b \end{array} \]^{p,L}_{a,b}  = \nn
&&  \widetilde{\sum}^{(0)} q^{\Phi^{p,0}_{a,b} ({\bf m}) + n_0 + 1-m_{b-1} +
m_b} \[ \begin{array}{l} {\bf n} + {\bf m} - {\bf E}_{0,1} - 2 {\bf
E}_{2,b-1} + {\bf e}_{b-1} - {\bf e}_b \\ {\bf n} +
{\bf e}_{b-1} - {\bf e}_b \end{array} \]^{p,L}_{a,b}  \nn
&& + \widetilde{\sum}^{(0)} q^{\Phi^{p,0}_{a,b} ({\bf m}) + n_0 + 1-m_{b-1} +
m_b + m_1} \[ \begin{array}{l} {\bf n} + {\bf m} - {\bf E}_{0,1} - 2 {\bf
E}_{2,b-1} + {\bf e}_{b-1} - {\bf e}_b \\ {\bf
n} - {\bf e}_1 + {\bf e}_{b-1} - {\bf e}_b \end{array} \]^{p,L}_{a,b}  \; .
\nn &&
\label{fourfortytwo}
\ena
Recalling (4.39), we identify the first place as $q^{L-1} F^{p,0}_{a,b}
(L-2)$.  Next we perform the change of
summation variables
\bea
 {\bf n} \rightarrow {\bf n} - {\bf e}_0 + {\bf e}_1 -{\bf e}_2 - {\bf
e}_{b-1} + {\bf e}_b  \; ;
 {\bf m} \rightarrow {\bf m} + 2 {\bf E}_{2,b-1}
\label{fourfortythree}
\ena
in the second piece to obtain
\bea
&& \widetilde{\sum}^{(0)} q^{\Phi^{p,0}_{a,b} ({\bf m}) + n_0 + 1-m_{b-1} +
m_b + m_1} \[ \begin{array}{l} {\bf n} + {\bf m} - {\bf E}_{0,1} - 2 {\bf
E}_{2,b-1} + {\bf e}_{b-1} - {\bf e}_b \\
{\bf n} - {\bf e}_1 + {\bf e}_{b-1} - {\bf e}_b \end{array} \]^{p,L}_{a,b}
= \nn
&& \widetilde{\sum}^{(0)} q^{\Phi^{p,0}_{a,b} ({\bf m}) + n_0 + m_2} \[
\begin{array}{l} {\bf n} + {\bf m} -
2 {\bf e}_0 - {\bf e}_2 \\ {\bf n} -  {\bf e}_0 - {\bf e}_2 \end{array}
\]^{p,L}_{a,b}
\label{fourfortyfour}
\ena
Hence by (4.39), (4.41), (4.42) and (4.44),
\bea
&& q^{L-\frac12} F^{p,0}_{a,b-1} (L-1) - q^{L-1} F^{p,0}_{a,b} (L-2) -
\widetilde{\sum}^{(0)} q^{\Phi^{p,0}_{a,b} ({\bf
m}) + n_0 + m_2} \[ \begin{array}{l} {\bf n} + {\bf m} - 2 {\bf e}_0  -
{\bf e}_2 \\ {\bf n} -  {\bf e}_0 - {\bf
e}_2 \end{array} \]^{p,L}_{a,b} = \nn
&&  \sum^{b-1}_{l=2}  \widetilde{\sum}^{(0)} q^{\Phi^{p,0}_{a,b} ({\bf m})
+ n_0 + 1-m_{b-1} +
m_b + (m_l -1)} \[ \begin{array}{l} {\bf n} + {\bf m} - {\bf e}_0  - {\bf
E}_{2,b-1} + {\bf e}_{b-1} - {\bf e}_b - {\bf
E}_{l,b-1} \\ {\bf n} +  {\bf e}_{b-1} - {\bf e}_b - {\bf e}_l \end{array}
\]^{p,L}_{a,b} \; . \nn &&
\label{fourfortyfive}
\ena
Let us now perform the $l$ - dependent change of summation variables
\bea
&& {\bf n} \rightarrow {\bf n} - {\bf e}_{l-1} +  {\bf e}_l + {\bf e}_{b-1}
- {\bf e}_b \nn
&& {\bf m} \rightarrow {\bf m} - 2 {\bf E}_{l,b-1}
\label{fourfortysix}
\ena
in the double sum featured in equation (4.35) to derive
\bea
&& \sum^b_{l=3}  \widetilde{\sum}^{(0)} q^{\Phi^{p,0}_{a,b} ({\bf m}) + n_0
+ m_l} \[ \begin{array}{l} {\bf n} +
{\bf m} - {\bf e}_0  - {\bf E}_{2,l}  \\ {\bf n} -  {\bf e}_l \end{array}
\]^{p,L}_{a,b} = \nn
&& \sum^{= b-1}_{l=2}  \widetilde{\sum}^{(0)} q^{\Phi^{p,0}_{a,b} ({\bf m})
+ n_0 + 1-m_{b-1} +
m_b + (m_l -1)} \[ \begin{array}{l} {\bf n} + {\bf m} - {\bf e}_0  - {\bf
E}_{2,b-1} + {\bf e}_{b-1} - {\bf e}_b - {\bf
E}_{l,b-1} \\ {\bf n} +  {\bf e}_{b-1} - {\bf e}_b - {\bf e}_l \end{array}
\]^{p,L}_{a,b} \; . \nn &&
\label{fourfortyseven}
\ena
We now use (4.36-37), (4.45) and (4.47) to simplify (4.35) as follows
\bea
F^{p,0}_{a,b} (L) &=& q^{L-1}  F^{p,0}_{a,b+1} (L-1) +  F^{p,0}_{a,b} (L-1)
+ q^{L-\frac12} F^{p,0}_{a,b-1}
(L-1)\nn
&-& q^{L-1}  F^{p,0}_{a,b} (L-2) - \widetilde{\sum}^{(0)}
q^{\Phi^{p,0}_{a,b} ({\bf m}) + n_0 + m_2} \[
\begin{array}{l} {\bf n} + {\bf m} - 2 {\bf e}_0 - {\bf e}_2 \\ {\bf n} -
{\bf e}_0 - {\bf e}_2 \end{array}
\]^{p,L}_{a,b} \nn
&+& \widetilde{\sum}^{(0)} q^{\Phi^{p,0}_{a,b} ({\bf m}) + n_0 + m_2} \[
\begin{array}{l} {\bf n} + {\bf
m} - {\bf e}_0 - {\bf e}_2\\ {\bf n}-  {\bf e}_2 \end{array} \]^{p,L}_{a,b}
\; .
\label{fourfortyeight}
\ena
The last two terms in the right-hand-side of (4.48) can be combined together
with the help of (2.30) as
\bea
\widetilde{\sum}^{(0)} q^{\Phi^{p,0}_{a,b} ({\bf m}) + 2 n_0 + m_2} \[
\begin{array}{l} {\bf n} + {\bf
m} - 2 {\bf e}_0  - {\bf e}_2\\ {\bf n}-  {\bf e}_2 \end{array} \]^{p,L}_{a,b}
\label{fourfortynine}
\ena
which according to (4.40) is nothing else but $q^{2(L-1)} F^{p,0}_{a,b}
(L-2)$.  Hence,
\bea
F^{p,0}_{a,b} (L) &=& q^{L-1}  F^{p,0}_{a,b+1} (L-1) +  F^{p,0}_{a,b} (L-1)
+ q^{L-\frac12} F^{p,0}_{a,b-1}
(L-1) \nn
&+& q^{L-1} (q^{L-1} -1) F^{p,0}_{a,b} (L-2) \; .
\label{fourfifty}
\ena
Thus, (4.5) is established.

\bigskip

\centerline{\underline{\large 4.4 Proof of even identities (2.12-13)}}
\bigskip
Comparing (3.17)-(3.22) and (4.1)-(4.6) one readily checks that polynomials
in the left-hand-side and right-hand-side
of (2.12-13) obey exactly the same recurrences, which specifies these
polynomials uniquely.  Therefore, to complete the
proof one needs to show that identities (2.12-13) hold for $L=0$.

Keeping in mind that all $n$ and $m$ variables in (2.19) are non-negative
and taking advantage of (2.25) and (2.28-29)
one concludes after a bit of analysis that the ${\bf n} - {\bf m}$ system
(2.21) with $L=0$ has solutions if and only if
$a=b$.  In the latter case the solution is
\bea
{\bf n} = {\bf e}_a ; {\bf m} =0 \; {\rm for} \; L=0, a=b
\label{fourfiftyone}
\ena
Substituting (4.51) into (2.19), yields
\bea
F^{p,0}_{a,b} (0) = \delta_{a,b} \; .
\label{fourfiftytwo}
\ena
On the other hand, using
\bea
\( \begin{array}{c} 0,q \\ A \end{array} \)^n = \delta_{A,0}
\label{fourfiftythree}
\ena
one immediately gets
\bea
B_{a,b} (0,s) = \tilde{B}_{a,b} (0,s) = \delta_{a,b} \; .
\label{fourfiftyfour}
\ena
From (2.10) one has
\bea
\varphi (a,a) = 0 \; .
\label{fourfiftyfive}
\ena
Hence by (4.52), (4.54) and (4.55), it is clear that (2.12-13) hold for
$L=0$.  Therefore, we completed the proof of
even identities (2.12-13).

\bigskip

\centerline{\underline{\large 4.5 Proof of odd identities (2.14-15)}}
\bigskip
Close examination of the proof of fermionic recurrences (4.1)-(4.6) given
above shows that the parity of the $m_{p-2}$ -
variable played a minor role in it.  One can literally repeat all the steps
in the
proof at (4.1)-(4.6) with parity of
$m_{p-2}$ reversed in all the earlier formulas to find ``odd" analogues
of the even recurrences. They are listed as follows:
\bea
\bar{F}^{p,1}_{p-a,1} (L) = \bar{F}^{p,1}_{p-a,1} (L-1) + q^{L-\frac12}
\bar{F}^{p,1}_{p-a,2} (L-1) \; {\rm for} \; a
\equiv 1 \; ({\rm mod} \; 2)
\label{fourfiftysix}
\ena
\bea
\bar{F}^{p,1}_{p-a,1} (L) = \bar{F}^{p,1}_{p-a,1} (L-1) + q^{L-1}
\bar{F}^{p,1}_{p-a,2} (L-1) \; {\rm for} \; a
\equiv 0 \; ({\rm mod} \; 2)
\label{fourfiftyseven}
\ena
\bea
\bar{F}^{p,1}_{p-a,2} (L) = q^{L-1} \bar{F}^{p,1}_{p-a,3} (L-1) +
q^{L-\frac12} \bar{F}^{p,1}_{p-a,1} (L-1) +
\bar{F}^{p,1}_{p-a,2} (L-1) \nn
{\rm for} \; a \equiv 1 \; ({\rm mod} \; 2) \hphantom{xxxxxxxxx}
\label{fourfiftyeight}
\ena
\bea
\bar{F}^{p,1}_{p-a,2} (L) = {1 \over \sqrt{q}} \;  \bar{F}^{p,1}_{p-a,3}
(L-1) +  \bar{F}^{p,1}_{p-a,1} (L-1) +
q^{L-1} \bar{F}^{p,1}_{p-a,2} (L-1) \nn
{\rm for} \; a \equiv 0 \; ({\rm mod} \; 2) \hphantom{xxxxxxxxx}
\label{fourfiftynine}
\ena
\bea
\bar{F}^{p,1}_{p-a,b} (L) = q^{L-1}  \bar{F}^{p,1}_{p-a,b+1} (L-1) +
q^{L-\frac12} \bar{F}^{p,1}_{p-a,b-1} (L-1) +
\bar{F}^{p,1}_{p-a,b} (L-1) \nn
+ q^{L-1} (q^{L-1} -1) \bar{F}^{p,1}_{p-a,b} (L-2) \hphantom{xxxxxxxxxxxxx} \nn
{\rm for} \; a+b \equiv 1 \; ({\rm mod} \; 2) , 2<b \leq p-1 \hphantom{xxx}
\label{foursixty}
\ena
\bea
\bar{F}^{p,1}_{p-a,b} (L) &=& {1 \over \sqrt{q}} \left\{
\bar{F}^{p,1}_{p-a,b+1} (L-1) + (q^{L-1}-1)
\bar{F}^{p,1}_{p-a,b+1} (L-2) \right\} \nn
&+& \left\{ \bar{F}^{p,1}_{p-a,b-1} (L-1) + (q^{L-1}-1)
\bar{F}^{p,1}_{p-a,b-1} (L-2) \right\} \nn
&+& \left\{ \bar{F}^{p,1}_{p-a,b} (L-1) + q^{L-2} (q^{L-1}-1)
\bar{F}^{p,1}_{p-a,b} (L-2) \right\} \nn
&& {\rm for} \; a+b \equiv 0 \; ({\rm mod} \; 2) , 2<b \leq p-1 \hphantom{xxx}
\label{foursixtyone}
\ena
with $2 \leq a \leq p-2, \bar{F}^{p,1}_{p-a,p} (L) \equiv 0$.  For $2 \leq
a \leq p-2, 1 \leq b \leq p-1$ polynomials
$\bar{F}^{p,1}_{p-a,b} (L)$ were defined as
\bea
\bar{F}^{p,1}_{p-a,b} (L) = \widetilde{\sum}^{(1)} q^{\Phi^{p,1}_{p-a,b} ({\bf
m})} \prod^{p-2}_{j=0} \[ \begin{array}{c} n_j
+ m_j \\ n_j \end{array} \]'_q
\label{foursixtytwo}
\ena
with ${\bf n}, {\bf m}$ variables being constrained by (2.21),(2.22) with
$a$ replaced by $p-a$ and $i=1$.  We would like to
remind the reader that as in the discussion given at the end of Section 2
\bea
\bar{F}^{p,1}_{p-a,b} (L) = F^{p,1}_{p-a,b} (L) + \delta_{L,0} \cdot
\delta_{a,b} \cdot \theta (1<b<p-1) \; .
\label{foursixtythree}
\ena
Once again, comparing (4.56)-(4.61) with (3.17)-(3.22) shows that
bosonic and fermionic polynomials in (2.14),(2.15)
obey exactly the same recursion relations.  Thus, to complete the proof it
is sufficient to demonstrate that (2.14)
holds for all $L \geq 0$.  This is easily done as follows.  First, from
definition (2.7), it is obvious that
\bea
F^{p,1}_{p-a,1} (L) = F^{p,0}_{p-a,p-1} (L) \; .
\label{foursixtyfour}
\ena
Next, replacing $a$ by $p-a$ and setting $b=p-1$ in the identity
(2.13), we get
\bea
F^{p,1}_{p-a,1} (L) = F^{p,0}_{p-a,p-1} (L) = \left\{ \begin{array}{l}
q^{\varphi (p-a, p-1)} B^p_{p-a,p-1} (L,p) \;
{\rm for} \; a \equiv 1 \; ({\rm mod} \; 2) \\ \\
q^{\tilde{\varphi} (p-a,p-1)} \tilde{B}^p_{p-a,p-1} (L,p+1) \; {\rm for} \; a
\equiv 0 \; ({\rm mod} \; 2) \end{array}
\right.
\label{foursixtyfive}
\ena
with $2 \leq a \leq p-2$.  Finally, making use of reflection properties
(3.23),(3.24) we derive (2.14).  Thus the odd
identities (2.14) and (2.15) are established for $2 \leq a \leq p-2, 1 \leq
b \leq p-1$.  Observing that
\bea
F^{p,1}_{p-1,b} (L) \equiv F^{p,0}_{1,b} (L) \; {\rm for} \; 1 \leq b \leq p-1
\label{foursixtysix}
\ena
and
\bea
\varphi (p-1,p-1) = \tilde{\varphi} (p-1, p-1) = 0 \; ,
\label{foursixtyseven}
\ena
we conclude that odd identities with $a=1, 1 \leq b \leq p-1$ reduce to
the even ones proven earlier in Section 4.4.

\mysection{Infinite L. Virasoro character identities}

In this section we will let $L \rightarrow \infty$ in the formulas
(2.12)-(2.15).  To this end it's convenient to
define fermionic objects
\bea
Fer^{p,i}_{a,b} (q) = \lim_{L \rightarrow \infty} F^{p,i}_{a,b} (L)
\label{fiveone}
\ena
with $1 \leq a \leq p-2; 1 \leq b \leq p-1; i=0,1$.  Employing the well-known
formula describing limiting behavior of $q$ -
binomials
\bea
\lim_{L \rightarrow \infty} \[ \begin{array}{c} L \\ k \end{array} \] = {1
\over (q)_k} \; {\rm with} \; L, k \in \zz_{\geq 0} \; ,
\label{fivetwo}
\ena
it is easy to derive explicit expressions for $Fer^{p,i}_{a,b} (q)$
\bea
Fer^{p,i}_{a,b} (q) = \sum_{\scriptstyle{\bf m}\atop\scriptstyle m_{p-2}
\equiv d(i) ({\rm mod} \: 2)} q^{{\frac{1}{4}} {\bf m} \tilde{C}_{p-1} {\bf
m} - {1 \over 2} {\bf A}^i_{a,b} \cdot {\bf m}}  {1 \over (q)_{m_0}}
\prod^{p-2}_{j=1} \[ \begin{array}{c} ( {1 \over
2} \tilde{I}_{p-1} {\bf m} + {\bf u}^i_{a,b} )_j \\ m_j \end{array} \]_q
\label{fivethree}
\ena
where
\bea
d(i) = \left\{ \begin{array}{c} i \; {\rm for} \; b \not= p-1 \\ 1-i \;
{\rm for} \; b=p -1 \end{array} \right.
\label{fivefour}
\ena
and the rest of notation is the same as in Section 1 and Section 2.
Equipped with definitions (5.1),(5.3) and using
limiting properties (3.25),(3.26) we can obtain character and extended
character identities.

\bigskip

\centerline{\underline{\large  Character identities}}
\bea
Fer^{p,0}_{a,1} (q) = \left\{ \begin{array}{l} q^{{a(a-1) \over 4}}
\chi^{p,p+1}_{a,1} (q) \; {\rm for} \; a +
\delta_{a,1} \equiv 0 \; ({\rm mod} \; 2) \\ q^{{(a-1)(a-2) \over 4}}
\chi^{p,p+1}_{a,2} (q) \; {\rm for} \; a +
\delta_{a,1} \equiv 1 \; ({\rm mod} \; 2) \end{array} \right.
\label{fivefive}
\ena
\bea
Fer^{p,1}_{p-a,1} (q) = \left\{ \begin{array}{l} q^{\varphi (p-a,p-1)}
\chi^{p,p+1}_{a,1} (q) \; {\rm for} \; a \equiv
1 \; ({\rm mod} \; 2) \\ q^{\tilde{\varphi}(p-a,p-1)} \chi^{p,p+1}_{a,2}
(q) \; {\rm for} \; a \equiv 0 \; ({\rm mod} \;
2)
\end{array} \right.
\label{fivesix}
\ena
\bea
Fer^{p,0}_{a,b} (q) = q^{\varphi(a,b)} \chi^{p,p+1}_{a,b+1} (q) \; {\rm
for} \; b+a+ \delta_{a,1}  \equiv
0 \; ({\rm mod} \; 2)
\label{fiveseven}
\ena
\bea
Fer^{p,1}_{p-a,b} (q) = q^{\tilde{\varphi}(p-a,p-1) + \varphi (a,b)-a {a-1
\over 4}} \chi^{p,p+1}_{a,b+1} (q) \; {\rm for}
\; b+a \equiv 0 \; ({\rm mod} \; 2)
\label{fiveeight}
\ena
with $1 \leq a \leq p-2 \; , \; 1 \leq b \leq p-1$\\
and
\bea
Fer^{p,1}_{p-1,b} (q) \equiv Fer^{p,0}_{1,b} (q) \; {\rm for} \; 1 \leq b
\leq p-1 \; .
\label{fivenine}
\ena

\bigskip

\centerline{\underline{\large Extended characters identities}}
\bea
Fer^{p,0}_{a,b} (q) = q^{\tilde{\varphi}(a,b)} \left\{ \theta (b \not= p-1)
\cdot \chi^{p,p+1}_{a,b+2} (q) \right.
\hphantom{xxxxx} \nn + \left. q^{a-b} \chi^{p,p+1}_{a,b} (q) \right\}
\hphantom{xxxxx} \nn
{\rm for} \; b+a + \delta_{a,1} \equiv 1 \; ({\rm mod} \; 2) \; , 2 \leq b
\leq p-1 \; ,\nn
1 \leq a \leq p-2 \hphantom{xxxxxxxxx}
\label{fiveten}
\ena
\bea
Fer^{p,1}_{p-a,b} (q) &=& q^{\tilde{\varphi}(p-a,p-1) + \tilde{\varphi} (a,b) -
{a -1 \over 4} (a+2)} \nn
&& \left\{ \theta (b \not= p-1) \chi^{p,p+1}_{a,b+2} (q) + q^{a-b}
\chi^{p,p+1}_{a,b} (q) \right\} \nn
&& \quad {\rm for} \; a+b \equiv 0 \; ({\rm mod} \; 2), \;  2 \leq b \leq
p-1 \nn
&& \quad {\rm and} \; 1 \leq a \leq p-2
\label{fiveten}
\ena
with $\chi^{p,p+1}_{a,b} (q)$ defined by (1.7).

S.O. Warnaar pointed out to us that the variable $m_1$ in (2.1) and
(2.7) can be summed out by means of the q-analogue of the
Chu-Vandermonde formula [(3.3.10) of \cite{and}]. 
Having done so and taking the limits $L\rightarrow
\infty$ we obtain from (2.13) with $2\leq b\leq p-2,~~ 2\leq a\leq p-2$
and $a+b\equiv 1({\rm mod}2)$ the following character identity
\bea
&&\sum_{{\bf m}\atop m_{p-2}\equiv 0 ({\rm mod}2)}
q^{{{\bf m}C_{p-2}{\bf m}\over 4}-{ m_1\over 2}-{ m_b\over 2}}{1\over
(q)_{m_1}}\prod_{j=2}^{p-2}{({1\over 2}{I}_{p-2}{\bf m})_j+{1\over
2}(\delta_{a,j}+\delta_{b,j})\atopwithdelims[] m_j}_q\nonumber \\
&&=q^{{\tilde\varphi}(a,b)}
\{\chi_{a,b+2}^{p,p+1}(q)+q^{a-b}\chi_{a,b}^{p,p+1}(q)\}
\ena
which was not noticed previously.

\bigskip

\mysection{Closing remarks}}

Here we would like to point out some salient features of the identities
(2.12)-(2.15) proven in this paper.  First, the
need to take into account certain negative solutions to the ${\bf n}, {\bf
m}$ system (2.21) as discussed at the end of
Section 2 makes our identities somewhat similar to those for non-unitary
minimal models $M (p, p^\prime \not= p+1)$
perturbed by the operator $\varphi_{1,3}$ \cite{btwentythree},
\cite{btwentythreea}.  Second, the bosonic forms
of our
polynomials required two types of $q$ -
trinomial coefficients.  The analogous feature was previously seen in the
study of $N=1$ superconformal models
$SM(p,p^\prime)$ where two types of $q$ - trinomials were used to describe
Neveu-Schwarz and Ramond sector identities
\cite{beighteen},\cite{bnineteen},\cite{btwentytwo},\cite{btwentyfour}.
This similarity is not very surprising because both  dilute and
$N=1$ superconformal models share common
spin-1 properties.

It is rather intriguing to notice, that despite the fact that we have a
complete set of Virasoro character identities
(5.5)-(5.8), we find in addition extended character identities
(5.9),(5.10).  In this regard, we want to mention that
linear combinations of Virasoro characters in (5.10),(5.11) were interpreted
by K. Seaton and L. Scott in \cite{bsixteen} as arising from
configuration sums for excited states.  Thus, it is conceivable that
extended character identities provide a hint for the
existence of some superselection rules which break up the space of our
models into a number of orthogonal subsectors.

Identities of this paper are associated with $A$ - modular invariant models
of conformal field theory in the classification
of \cite{btwentyfive}.  The only  known result 
for $D$ - modular invariant models
perturbed by operator $\varphi_{2,1}$ is for
three-state Potts model $M_D$(5,6).  In this case the identities are of a
parafermionic form $Z_3$ introduced by J. Lepowsky
and M. Primc \cite{btwentysix} and are quite different 
from the ones studied here.  In
particular, the para-fermionic form has two particles
instead of the four particles found here for $M$(5,6).  Thus, one may
speculate that not only different integrable
perturbations lead to distinct Rogers-Ramanujan type identities but also
different modular invariance specifications
may yield distinct new identities.  Finally, it is highly desirable to find a
Partition Theoretical interpretation of
polynomials $B^p_{a,b} (L,s), \tilde{B}^p_{a,b} (L,s)$.  Such an
interpretation may lead to the $q$ - trinomial
generalization of the Burge iterative process \cite{btwentyseven} 
which was recently used
by O. Foda, et al. \cite{btwentyeight} to give another proof
of Rogers-Ramanujan type identities for $M(p, p^\prime)$ models perturbed by
operator $\varphi_{1,3}$.

\bigskip

{\Large \bf {Note added}} After the completion of this paper
S.O.Warnaar informed us that he has obtained an independent proof of
identity (\ref{oneten}).

\bigskip

{\large \bf{Acknowledgment}}  This work is supported in part by the
National Science Foundation grant of the USA under DMR9703543.We wish
to thank S.O. Warnaar for fruitful discussions.

\appendix\myappendix

Here we will show that under replacement
\bea
\[ \begin{array}{l} n_0 + m_0 \\ n_0 \end{array} \]_q \rightarrow \[
\begin{array}{l} n_0 + m_0 \\ n_0 \end{array} \]^\prime_q
\ena
in (2.19), the combined contribution due to the negative values of $n_0$
vanishes,
unless $L=0; p-a=b \not= 1, p-1; i=1$ in
(2.19).  Since the main features of our treatment remain essentially the
same for all $1 \leq a \leq p-2; 1 \leq b \leq
p-1; i=0,1$ in (2.19), we first limit our attention to the case $a=b=1;
i=0$ in order to simplify a bit the discussion
below.  Note, that in this case equation (2.27) becomes
\bea
L= (n_0 + n_1) + m_2 +2 \sum^{p-2}_{j=2} n_j + m_{p-2} \; .
\ena
For $n_0$ to be negative
\bea
m_2 > 2L
\ena
should hold.  Evaluating (2.19) with replacement (A.1) for $m_2 > 2L$, one
encounters the sum
\bea
\sum^{C_0}_{m_1 =0} (-)^{m_1}  q^{\frac{m^2_1 - m_1}{2} - Lm_1} \[
\begin{array}{l} {m_2 \over 2} -L-1 + m_0 \\ {m_2
\over 2} - L-1 \end{array} \]_q  \[ \begin{array}{l} {m_2 \over 2} \\ m_1
\end{array} \]_q
\ena
with variables $m_0, m_1$ being constrained by $m_0 + m_1 =C_0$.  Below we
will show that this sum becomes zero, unless
\bea
0 \leq C_0 = m_0 + m_1 \leq L \; .
\ena
However, condition (A.5) is incompatible with (A.2), (2.34), (2.35).
Indeed, from (A.5) and
(2.26) it is readily seen that
\bea
n_0 + n_1 \geq 0 \; .
\ena
Now combining (A.2), (A.3), (A.6), (2.34) and (2.35) we arrive at contradiction
\bea
L > 2 L
\ena
which proves our claim.

It remains to demonstrate that sum (A.4) will vanish for $C_0 > L$.  To
this end let
\bea
P(L,C_0, B) = \sum^{\infty}_{m_1=0} (-)^{m_1} q^{\frac{m^2_1 - m_1}{2} -
Lm_1} \[ \begin{array}{l} B+C_0 - m_1 \\ B
\end{array} \]_q \[ \begin{array}{c} B+1+L \\ m_1 \end{array} \]_q \; .
\ena

Using the q-Chu-Vandermonde sum (equation II.6 in \cite{btwentyone})
the expression (A.8) is evaluated as
\bea
P(L,C_0,B)=(-1)^{C_0}q^{{C_0(C_0-1)\over 2}-C_0 L}{L\atopwithdelims[] C_0}_q.
\ena
Clearly $P(L,C_0,B)$ vanishes for $C_0>L.$
Finally, observing that (A.4) is nothing else but $P(L, C_0, \frac{m_2}{2}
- L-1)$, we are done.

Proceeding as above in general case: $ 1 \leq a \leq p-2; 1 \leq b \leq
p-1; i=0, 1$, we conclude that contributions due
to negative values of $n_0$ must vanish unless (A.6) holds and
\bea
m_2 \geq \left\{ \begin{array}{l} 2L +1 \; {\rm if} \; m_2 \; {\rm is \; odd}
\\ 2L+2 \; {\rm if} \; m_2 \; {\rm is \; even}
\end{array} \right. \; .
\ena
Again, combining (2.27), (A.6), (A.10), (2.34) and (2.35) we arrive at
a contradiction, unless $L \leq 1$.  Detailed
analysis of the remaining cases with $L=0,1$ shows that one can have
an overall contribution due to negative values of
$n_0$ if and only if $L=0; p-a=b \not= 1, p-1$ and $i=1$ in (2.19).  In the
latter case one needs to take into account
the following solution to (2.21) with negative $n_0$:
\bea
n_0 = - n_1 =-1
\ena
\bea
n_i = 0 \; {\rm for} \; 2 \leq i \leq p-2
\ena
\bea
m_0 = m_1 =0
\ena
and for $i = 2,3, \cdots, p-2$
\bea
m_i = \left\{ \begin{array}{l}  (i-1) \; {\rm for} \; 2 \leq i \leq \; {\rm
min} \;  \{ b, p-b \} \\(p-1-i) \; {\rm for
\; max} \; \{ b,p-b \} \leq i \leq p-2 \\ {\rm min} \; \{ p-b, b \} -1 \;
{\rm otherwise} \; . \end{array} \right.
\ena
It can be easily checked that this solution contributes as 1.

\myappendix

Here we will prove recurrences (4.6).  Throughout this appendix it is
assumed that $b+a+ \delta_{a,1} \equiv 1$ (mod 2)
and $1 \leq a \leq p-2, 2< b \leq p-1$.  As usual, we start with telescopic
expansion for $F^{p,0}_{a,b} (L)$.
\bea
F^{p,0}_{a,b} (L) &=& \theta (b \not= p-1) \widetilde{\sum}^{(0)}
q^{\Phi^{p,0}_{a,b} ({\bf m})}
\[ \begin{array}{l} {\bf n} + {\bf m} - {\bf E}_{1,b} \\
{\bf n}  \end{array} \]^{p,L}_{a,b}  + \nn
&+& \sum^b_{l=3} \widetilde{\sum}^{(0)} q^{\Phi^{p,0}_{a,b} ({\bf m}) + m_l}
 \[ \begin{array}{l} {\bf n} + {\bf m} - {\bf E}_{1,l}  \\ {\bf n} -
{\bf e}_l  \end{array} \]^{p,L}_{a,b}  \nn
&+& \sum^2_{l=1} \widetilde{\sum}^{(0)} q^{\Phi^{p,0}_{a,b} ({\bf m}) + m_l}
\[ \begin{array}{l} {\bf n} + {\bf m} - {\bf E}_{1,l} \\ {\bf
n} - {\bf e}_l  \end{array} \]^{p,L}_{a,b}
\ena
with $m_{p-1} \equiv 0$.  Next we perform the $l$ - dependent change of
summation variables
\bea
{\bf n} \rightarrow {\bf n} - {\bf e}_{l-1} + {\bf e}_l + {\bf e}_{b-1} -
{\bf e}_b ; \; {\bf m} \rightarrow {\bf m} - 2 {\bf
E}_{l, b-1}
\ena
in the first double sum in (B.1) to find after replacing $l-1$ by $l$ that
\bea
&& \sum^b_{l=3} \widetilde{\sum}^{(0)} q^{\Phi^{p,0}_{a,b} ({\bf m}) + m_l}
 \[ \begin{array}{l} {\bf n} + {\bf m} - {\bf E}_{1,l}  \\ {\bf n} -
{\bf e}_l  \end{array} \]^{p,L}_{a,b}  = \nn
&& = \sum^{b-1}_{l=2} \widetilde{\sum}^{(0)} q^{\Phi^{p,0}_{a,b} ({\bf m}) +
m_b - m_{b-1} +m_l}
 \[ \begin{array}{l} {\bf n} + {\bf m} - {\bf E}_{1,b-1} -  {\bf E}_{l,b-1}
+ {\bf e}_{b-1} - {\bf e}_b \\ {\bf n} -
{\bf e}_l + {\bf e}_{b-1} - {\bf e}_b \end{array} \]^{p,L}_{a,b}  \; . \nn &&
\ena
To proceed further we will use
\bea
F^{p,0}_{a,b+1} (L-1) &=& \sqrt{q} \; \widetilde{\sum}^{(0)}
q^{\Phi^{p,0}_{a,b} ({\bf m})}
\[ \begin{array}{l} {\bf n} + {\bf m} - {\bf E}_{1,b}  \\
{\bf n}  \end{array} \]^{p,L}_{a,b} \; , \\
F^{p,0}_{a,b-1} (L-1) &=& \widetilde{\sum}^{(0)} q^{\Phi^{p,0}_{a,b} ({\bf m})
+ m_b - m_{b-1} +1}
\[ \begin{array}{l} {\bf n} + {\bf m} - {\bf E}_{1,b-1} + {\bf e}_{b-1} -
{\bf e}_b  \\
{\bf n} + {\bf e}_{b-1} - {\bf e}_b \end{array} \]^{p,L}_{a,b}  \; ,  \\
q^{L-2} F^{p,0}_{a,b} (L-2) &=& \widetilde{\sum}^{(0)} q^{\Phi^{p,0}_{a,b}
({\bf m}) + m_b - m_{b-1} +1 +n_0}
\[ \begin{array}{l} {\bf n} + {\bf m} - {\bf E}_{0,b-1} - {\bf E}_{2,b-1} +
{\bf e}_{b-1} - {\bf e}_b \\
{\bf n} + {\bf e}_{b-1} - {\bf e}_b \end{array} \]^{p,L}_{a,b} \; , \nn && \\
q^{L-1} F^{p,0}_{a,b} (L-1) &=& \widetilde{\sum}^{(0)} q^{\Phi^{p,0}_{a,b}
({\bf m}) + n_0 + m_1}
\[ \begin{array}{l} {\bf n} + {\bf m} - {\bf e}_0 \\ {\bf n} - {\bf e}_1
\end{array} \]^{p,L}_{a,b} \; .
\ena
The telescopic expansion of (B.5) yields
\bea
&& F^{p,0}_{a,b-1} (L-1) =  \widetilde{\sum}^{(0)} q^{\Phi^{p,0}_{a,b} ({\bf
m}) + m_b -m_{b-1} +1}
\[ \begin{array}{l} {\bf n} + {\bf m} - {\bf e}_1 - 2  {\bf E}_{2,b-1} +
{\bf e}_{b-1} - {\bf e}_b\\
{\bf n} + {\bf e}_{b-1} - {\bf e}_b \end{array} \]^{p,L}_{a,b} + \nn
&&  \sum^{b-1}_{l=2} \widetilde{\sum}^{(0)} q^{\Phi^{p,0}_{a,b} ({\bf m}) +
m_b - m_{b-1} +m_l}
 \[ \begin{array}{l} {\bf n} + {\bf m} - {\bf E}_{1,b-1} -  {\bf E}_{l,b-1}
+ {\bf e}_{b-1} - {\bf e}_b \\ {\bf n} -
{\bf e}_l + {\bf e}_{b-1} - {\bf e}_b \end{array} \]^{p,L}_{a,b}  \; .
\ena
Combining (B.1), (B.4), (B.3) and (B.8) one readily derives
\bea
&& F^{p,0}_{a,b} (L) = {1 \over \sqrt{q}} F^{p,0}_{a,b+1} (L-1)  +
F^{p,0}_{a,b-1}
(L-1) \nn
&& \qquad - \widetilde{\sum}^{(0)} q^{\Phi^{p,0}_{a,b} ({\bf m}) + m_b -
m_{b-1} +1}
 \[ \begin{array}{l} {\bf n} + {\bf m} - {\bf e}_1 - 2 {\bf E}_{2,b-1} +
{\bf e}_{b-1}
- {\bf e}_b \\ {\bf n} +
{\bf e}_{b-1} - {\bf e}_b \end{array} \]^{p,L}_{a,b}  \nn
&& \qquad + \sum^2_{l=1} \widetilde{\sum}^{(0)} q^{\Phi^{p,0}_{a,b} ({\bf m}) +
m_l}
 \[ \begin{array}{l} {\bf n} - {\bf E}_{1,l} \\ {\bf n} -
{\bf e}_l \end{array} \]^{p,L}_{a,b} \; .
\ena
If we now replace $L$ by $L-1$ in (B.9) and use
\bea
\{ {\bf n}, {\bf m} \}^{L-1,0}_{a,b} = \{ {\bf n}, {\bf m} \}^{L,0}_{a,b} -
\{ {\bf e}_0, 0 \} \; ,
\ena
we obtain
\bea
&& F^{p,0}_{a,b} (L-1) = {1 \over \sqrt{q}} \;  F^{p,0}_{a,b+1} (L-2)  +
F^{p,0}_{a,b-1} (L-2) \nn
&& \qquad - \widetilde{\sum}^{(0)} q^{\Phi^{p,0}_{a,b} ({\bf m}) + m_b -
m_{b-1} +1}
 \[ \begin{array}{l} {\bf n} + {\bf m} - {\bf E}_{0,1} -2 {\bf E}_{2,b-1} +
{\bf e}_{b-1} - {\bf e}_b \\ {\bf n} -
{\bf e}_0 + {\bf e}_{b-1} - {\bf e}_b \end{array} \]^{p,L}_{a,b}  \nn
&& \qquad + \sum^2_{l=1} \widetilde{\sum}^{(0)} q^{\Phi^{p,0}_{a,b} ({\bf m}) +
m_l}
 \[ \begin{array}{l} {\bf n} - {\bf e}_0 - {\bf E}_{1,l} \\ {\bf n} -
{\bf e}_0 - {\bf e}_l\end{array} \]^{p,L}_{a,b} \; .
\ena
Next subtract (B.11) from (B.9) and use (2.30) and (B.6) to get
\bea
F^{p,0}_{a,b} (L) -  F^{p,0}_{a,b} (L-1) &=& {1 \over \sqrt{q}} \; \left\{
F^{p,0}_{a,b+1} (L-1) - F^{p,0}_{a,b+1}
(L-2) \right\} + \nn
&+& \left\{ F^{p,0}_{a,b-1} (L-1) -  F^{p,0}_{a,b-1} (L-2)  \right\} - \nn
&-& q^{L-2} F^{p,0}_{a,b-1} (L-2) + \nn
&+& \sum^2_{l=1}  \widetilde{\sum}^{(0)} q^{\Phi^{p,0}_{a,b} ({\bf m}) + n_0
+ m_l}
 \[ \begin{array}{l} {\bf n} - {\bf E}_{0,l} \\ {\bf n} - {\bf e}_l
\end{array} \]^{p,L}_{a,b} \; .
\ena
It remains to process the double sum in (B.12).  To this end we expand
(B.7) with the help of (2.31) as
\bea
q^{L-1} F^{p,0}_{a,b} (L-1) &=& \widetilde{\sum}^{(0)} q^{\Phi^{p,0}_{a,b}
({\bf m}) + n_0 + m_1}
\[ \begin{array}{l} {\bf n} + {\bf m} - {\bf E}_{0,1} \\
{\bf n} - {\bf e}_1  \end{array} \]^{p,L}_{a,b} \nn
&+& \widetilde{\sum}^{(0)} q^{\Phi^{p,0}_{a,b} ({\bf m}) + n_0 + 2 m_1 +1}
 \[ \begin{array}{l} {\bf n}+ {\bf m} - {\bf E}_{0,1} \\ {\bf n} - 2 {\bf
e}_1 \end{array} \]^{p,L}_{a,b}
\ena
and perform a change of summation variables
\bea
{\bf n} \rightarrow  {\bf n} - {\bf e}_0 + {\bf e}_1 \; ; \; {\bf m}
\rightarrow {\bf m} + {\bf e}_0 - {\bf e}_1
\ena
in the last sum in (B.13) to find
\bea
&& \widetilde{\sum}^{(0)} q^{\Phi^{p,0}_{a,b} ({\bf m}) + n_0 + 2 m_1 +1}
 \[ \begin{array}{l} {\bf n}+ {\bf m} - {\bf E}_{0,1} \\ {\bf n} - 2 {\bf
e}_1 \end{array} \]^{p,L}_{a,b} = \nn
&& \quad q^{L-1}\widetilde{\sum}^{(0)} q^{\Phi^0_{a,b} ({\bf m}) + m_1}  \[
\begin{array}{l} {\bf n} + {\bf m} - {\bf
E}_{0,1} \\ {\bf n} - {\bf E}_{0,1}  \end{array} \]^{p,L}_{a,b} \; .
\ena
Making use of (B.11), (B.13) and (B.15), it is readily seen that
\bea
&& \widetilde{\sum}^{(0)} q^{\Phi^{p,0}_{a,b} ({\bf m}) + n_0 + m_1}
 \[ \begin{array}{l} {\bf n}+ {\bf m} - {\bf E}_{0,1} \\ {\bf n} -  {\bf
e}_1 \end{array} \]^{p,L}_{a,b} = \nn
&& q^{L-1} \left\{ {1 \over \sqrt{q}} \; F^{p,0}_{a,b+1} (L-2) +
F^{p,0}_{a,b-1} (L-2)  \right.  \nn
&& \quad  - \widetilde{\sum}^{(0)} q^{\Phi^{p,0}_{a,b} ({\bf m}) + m_b -
m_{b-1} + 1}
 \[ \begin{array}{l} {\bf n} + {\bf m} - {\bf E}_{0,1} -  2 {\bf E}_{2,b-1}
+ {\bf e}_{b-1} - {\bf e}_b \\ {\bf n} -
{\bf e}_0 + {\bf e}_{b-1} - {\bf e}_b \end{array} \]^{p,L}_{a,b}  \nn
&& \quad + \widetilde{\sum}^{(0)} q^{\Phi^{p,0}_{a,b} ({\bf m}) + m_2}
 \left. \[ \begin{array}{l} {\bf n} + {\bf m} - {\bf E}_{0,2}  \\ {\bf n} -
{\bf e}_0 - {\bf e}_2 \end{array}
\]^{p,L}_{a,b}  \right\}  \; .
\ena
The last equation together with (B.12) yields
\bea
F^{p,0}_{a,b} (L) &=& {1 \over \sqrt{q}} \left\{ F^{p,0}_{a,b+1} (L-1) +
(q^{L-1} -1) F^{p,0}_{a,b+1} (L-2)  \right\}
+\nn
&+& \left\{ F^{p,0}_{a,b-1} (L-1) +  (q^{L-1} -1) F^{p,0}_{a,b-1} (L-2)
\right\} \nn
&+& F^{p,0}_{a,b} (L-1) + q^{L-2} (q^{L-1} -1) F^{p,0}_{a,b} (L-2)  \nn
&+& Z (L,a,b)
\ena
where
\bea
Z (L,a,b) = - q^{L-1} X (L,a,b) + Y (L,a,b)
\ena
with
\bea
X(L,a,b) = q^{L-2} F^{p,0}_{a,b} (L-2) + \widetilde{\sum}^{(0)}
q^{\Phi^{p,0}_{a,b} ({\bf m}) + m_b - m_{b-1} + 1}
 \[ \begin{array}{l} {\bf n}+ {\bf m} - {\bf E}_{0,1} -2 {\bf E}_{2,b-1} +
{\bf e}_{b-1} - {\bf e}_b \\ {\bf n} -
{\bf e}_0 + {\bf e}_{b-1} - {\bf e}_b \end{array} \]^{p,L}_{a,b} \nn &&
\ena
and
\bea
Y(L,a,b) &=& \widetilde{\sum}^{(0)} q^{\Phi^{p,0}_{a,b} ({\bf m}) + n_0 + m_2}
\[ \begin{array}{l} {\bf n} + {\bf m} - {\bf
E}_{0,2} \\ {\bf n} - {\bf e}_2 \end{array} \]^{p,L}_{a,b}  \nn
&+& \widetilde{\sum}^{(0)} q^{\Phi^{p,0}_{a,b} ({\bf m}) + L-1 + m_2}
 \[ \begin{array}{l} {\bf n} + {\bf m} - {\bf E}_{0,2}  \\ {\bf n} - {\bf
e}_0 -{\bf e}_2 \end{array}
\]^{p,L}_{a,b} \; .
\ena
Making use of (B.6) and (2.30) we can simplify (B.19) to
\bea
X(L,a,b) =  \widetilde{\sum}^{(0)} q^{\Phi^{p,0}_{a,b} ({\bf m}) + m_b -
m_{b-1}
+ 1}
 \[ \begin{array}{l} {\bf n}+ {\bf m} - {\bf e}_1 -2 {\bf E}_{2,b-1} + {\bf
e}_{b-1} - {\bf e}_b \\ {\bf n} +
{\bf e}_{b-1} - {\bf e}_b \end{array} \]^{p,L}_{a,b} \; . \nn &&
\ena
Next we perform the change of summation variables
\bea
{\bf n} \rightarrow  {\bf n} + {\bf e}_b - {\bf e}_{b-1} -  {\bf e}_2 \; ;
\; {\bf m} \rightarrow {\bf m} + 2 {\bf
E}_{2,b-1} - {\bf e}_0 + {\bf e}_1
\ena
in (B.21) to get
\bea
X(L,a,b) =  \widetilde{\sum}^{(0)} q^{\Phi^{p,0}_{a,b} ({\bf m}) - m_0 +
{m_2 +1
\over 2}}
 \[ \begin{array}{l} {\bf n}+ {\bf m} - {\bf e}_0 - {\bf e}_2 \\ {\bf n} -
{\bf e}_2 \end{array} \]^{p,L}_{a,b} \; .
\ena
Finally, multiplying both sides in (B.23) by $q^{L-1}$ and using (4.10) gives
\bea
q^{L-1} X(L,a,b) = \widetilde{\sum}^{(0)} q^{\Phi^{p,0}_{a,b} ({\bf m}) +
n_0 + m_2}
 \[ \begin{array}{l} {\bf n}+ {\bf m} - {\bf e}_0 - {\bf e}_2 \\ {\bf n} -
{\bf e}_2 \end{array} \]^{p,L}_{a,b} \; .
\ena
It remains to process (B.20).  To this end we perform the change of
summation variables
\bea
{\bf n} \rightarrow  {\bf n} + {\bf e}_0 - {\bf e}_1 \; ; \; {\bf m}
\rightarrow {\bf m} -{\bf
e}_0 + {\bf e}_1
\ena
in the second sum in (B.20) to find
\bea
&& \widetilde{\sum}^{(0)} q^{\Phi^{p,0}_{a,b} ({\bf m}) + L-1 - m_2}
 \[ \begin{array}{l} {\bf n}+ {\bf m} - {\bf E}_{0,2} \\ {\bf n} -
{\bf e}_0 - {\bf e}_2 \end{array} \]^{p,L}_{a,b} = \nn
&& \quad \widetilde{\sum}^{(0)} q^{\Phi^{p,0}_{a,b} ({\bf m}) + n_0 + m_2 +
m_1}
 \[ \begin{array}{l} {\bf n}+ {\bf m} - {\bf E}_{0,2} \\ {\bf n} -
{\bf E}_{1,2} \end{array} \]^{p,L}_{a,b}  \; .
\ena
We can now combine the two terms in (B.20) with the help of (2.31) as
\bea
Y(L,a,b) =  \widetilde{\sum}^{(0)} q^{\Phi^{p,0}_{a,b} ({\bf m}) + n_0 + m_2}
 \[ \begin{array}{l} {\bf n}+ {\bf m} - {\bf e}_0 - {\bf e}_2 \\ {\bf n} -
{\bf e}_2 \end{array} \]^{p,L}_{a,b} \; .
\ena
The desired result (4.6) now easily follows from (B.17), (B.18), (B.24),
(B.27).  Thus, recurrences (4.6) are
established.

\bibliographystyle{plain}

\end{document}